# Automated Search for Conjectures on Mathematical Constants using Analysis of Integer Sequences

Ofir Razon, Yoav Harris, Shahar Gottlieb, Dan Carmon, Ofir David and Ido Kaminer

Technion - Israel Institute of Technology, Haifa 3200003, Israel

**Abstract**

Formulas involving fundamental mathematical constants had a great impact on various fields of science and mathematics, for example aiding in proofs of irrationality of constants. However, the discovery of such formulas has historically remained scarce, often perceived as an act of mathematical genius by great mathematicians such as Ramanujan, Euler, and Gauss. Recent efforts to automate the discovery of formulas for mathematical constants, such as the Ramanujan Machine project, relied on exhaustive search. Despite several successful discoveries, exhaustive search remains limited by the space of options that can be covered and by the need for vast amounts of computational resources. Here we propose a fundamentally different method to search for conjectures on mathematical constants: through analysis of integer sequences. We introduce the Enumerated Signed-continued-fraction Massey Approve (ESMA) algorithm, which builds on the Berlekamp-Massey algorithm to identify patterns in integer sequences that represent mathematical constants. The ESMA algorithm found various known formulas for $e, e^2, \tan(1)$, and ratios of values of Bessel functions. The algorithm further discovered a large number of new conjectures for these constants, some providing simpler representations and some providing faster numerical convergence than the corresponding simple continued fractions. Along with the algorithm, we present mathematical tools for manipulating continued fractions. Specifically, we present novel transformations of continued fractions with alternating positive and negative signs in the numerators, showing how to convert them to polynomial continued fractions and to simple continued fractions. These connections enable us to characterize what space of constants can be found by ESMA and quantify its algorithmic advantage in certain scenarios. Altogether, this work continues in the development of augmenting mathematical intuition by computer algorithms, to help reveal mathematical structures and accelerate mathematical research.

## 1 Introduction

Fundamental mathematical constants like $e, \pi,$ and the Golden Ratio $\varphi$ are ubiquitous in almost all fields of science and mathematics [1]. The discovery of new formulas involving fundamental mathematical constants often inspired mathematical research that revealed intrinsic properties of said constants and occasionally had a great impact in seemingly unrelated fields. It has been said that Ramanujan's 'Lost Notebook', which contains more than 600 conjectures (many involving fundamental constants) [2], has "caused roughly as much stir in the mathematical world as the discovery of Beethoven's Tenth symphony would cause in the musical world" [3]. The formulas in the 'Lost Notebook' about q-series and on mock modular forms were found to have countless applications in physics such as calculating black hole entropy [4]. The impact of formulas on fundamental mathematical constants is further exemplified with Apery's use of a polynomial continued fraction (CF) formula of $\zeta(3)$ to prove its irrationality [5]. As such, formulas involving mathematical constants often have implications beyond the formulas themselves, leading to the discovery of new mathematical structures.

Despite their great impact, the discovery of new formulas has been a scarce occasion, relying mostly on the cleverness and unfathomable intuition of great mathematicians such as Ramanujan. In the Ramanujan Machine Project [6], an automated approach for the discovery of conjectures involving fundamental constants was presented. The automation of conjecture discovery was approached as a search problem: different functions of given fundamental constant (e.g.,$1 + e/1 - e$) are calculated to some limited precision and stored in hash tables. Then, the algorithm calculates decimal values of a generated family of polynomial CFs (a CF in which the partial numerator and denominator sequences are integer-coefficient polynomials [7]) and searches for numerically matching values in the hash tables. Each match is further validated by a high decimal precision calculation. This algorithm obtained significant results, discovering many previously unknown formulas and conjectures for $e, \pi, \pi^2,$ and $\zeta(3)$, along with discovering the most efficient representation of the Catalan constant $G$ [6]. However, the algorithm was limited in its capacities since in its core, it still relied on an exhaustive search, expensive in computational resources and limited in the space of options that it can cover.

Here we propose a fundamentally different approach, based on the conversion of each target constant into a set of integer sequences, which we then analyze with an algorithm that attempts to detect a pattern in each of them. Such a pattern, if exists, may provide a formula for the target constant. The pattern recognition is based on a novel use of the **Berlekamp-Massey** algorithm [8-9]. We thus name the overall algorithm **Enumerated Signed-CF Massey Approve (ESMA)**.

Our approach is inspired by ideas of compression, entropy, and information theory [10]. Decimal representations of irrational constants such as e may seem (incorrectly) to contain an infinite amount of data with no discernable pattern, e.g., infinite entropy. Meaning, transmitting e (or any irrational constant) may seem to require a transmission of an infinite sequence of random digits to fully describe the constant. However, since there exists a formula to calculate e to infinite accuracy, its decimal representation actually contains zero entropy. In the example of e, this fact can be seen by analyzing its simple CF expansion [11], which shows a clear pattern $(\ldots 1, 2k, 1, 1, 2k+2, 1, 1, 2k+4, \ldots)$:

$$e - 2 = \cfrac{1}{1+\cfrac{1}{2+\cfrac{1}{1+\cfrac{1}{1+\cfrac{1}{4+\cfrac{1}{1+\cfrac{1}{1+\cfrac{1}{6+\cfrac{1}{1+\cfrac{1}{1+\cfrac{1}{8+\cfrac{1}{1+\cfrac{1}{\cdots}}}}}}}}}}}}}. \quad (1)$$

As such, to transmit e we can simply transmit the sequence it follows in its simple CF expansion. The vision behind our algorithm is to efficiently identify such patterns. We expand each target constant to a set of different CFs (with alternating positive and negative signs in the numerators) with the hope that at least one of them will reveal a pattern. Each successful result can be seen as a different compression of the seemingly infinite decimal representation of a constant to a zero-entropy formula.

We demonstrate the potential of this algorithm with equations involving **Signed Interlaced Continued Fractions (SICFs)** in which the partial denominator is a sequence made of $\beta$ different sub-sequences (interlaced sequence), and the partial numerator is some periodic sequence of $\pm 1$:

$$a_0 + \mathbb{K}_1^\infty \frac{b_j}{a_j} = a_0 + \cfrac{b_1}{a_1+\cfrac{b_2}{a_2+\cfrac{b_3}{a_3+\cdots}}} = a_0 + \cfrac{1}{a_1-\cfrac{1}{a_2-\cfrac{1}{a_3+\cdots}}}, \quad (2)$$

where $a_j \in \mathbb{Z}$ , $b_j \in \{1, -1\}$ $\forall j \in \mathbb{N}$, are the partial denominators and numerators of the CF, respectively. We specifically choose to utilize SICFs as their signed partial numerator $b_j = \pm 1$ can be used to directly extract an integer sequence given a fundamental constant. We attempt to recognize patterns in the obtained extracted integer sequence, $a_j$.

The ESMA algorithm was able to produce various known mathematical formulas, e.g.,

$$\tan(1) = 1 + \cfrac{1}{1+\cfrac{1}{1+\cfrac{1}{3+\cfrac{1}{1+\cfrac{1}{5+\cdots}}}}}, \; b_j = 1, \; a_j = \begin{cases} 1 & j = 2k+1 \\ 1+2k & j = 2k+2 \end{cases} \quad (3)$$

Along with a plethora of novel conjectures, some of which converge faster than their simple CFs e.g.,

$$\frac{2+2e}{-1+3e} = 2 - \cfrac{1}{1+\cfrac{1}{24+\cfrac{1}{3-\cfrac{1}{2+\cdots}}}}, \; b_j = \begin{cases} -1 & j = 6k+1 \\ 1 & j = 6k+2 \\ 1 & j = 6k+3 \\ -1 & j = 6k+4 \\ 1 & j = 6k+5 \\ 1 & j = 6k+6 \end{cases}, \; a_j = \begin{cases} 2 & j = 6k \\ 1+4k & j = 6k+1 \\ 24+64k & j = 6k+2 \\ 3+4k & j = 6k+3 \\ 2 & j = 6k+4 \\ 13+16k & j = 6k+5 \end{cases} \quad (4)$$

$$\frac{J_1(1)}{J_3(1)} = 23 - \cfrac{1}{1+\cfrac{1}{1+\cfrac{1}{39-\cfrac{1}{2+\cdots}}}}, \; b_j = \begin{cases} -1 & j = 3k+1 \\ 1 & j = 3k+2 \\ 1 & j = 3k+3 \end{cases}, \; a_j = \begin{cases} 23+16k & j = 3k \\ 1+k & j = 3k+1 \\ 1 & j = 3k+2 \end{cases}$$

$$\forall j \in \mathbb{N}, k \in \{0,1,2,3 \ldots\}$$

$J_y(x)$ are the Bessel functions of the first kind of order $y$ with argument $x$. Additional results are provided in Appendix A, demonstrating the potential of the algorithm.

Conjectures found are verified to a precision of 1000 decimal places. Therefore, obtaining a false conjecture that is merely a mathematical coincidence is highly unlikely. Specifically, every conjecture is found and verified in two steps: (1) Using the Berlekamp-Massey algorithm, a pattern is found in the sequence $a_j$ that we generated directly from the constant. We usually test 50 elements of the sequence $a_j$. More elements can be taken in cases where we want to find longer patterns. (2) If a pattern was found by the algorithm, additional validation is acquired by using the pattern to calculate additional elements of the $a_j$ sequence. We then evaluate the decimal representation of the CF and compare it with additional digits of the constant. We usually compare to the next 1000 digits of the constant's decimal representation. Consequently, the probability for a false positive is the probability that a random set of parameters satisfies an equivalence for 1000 digits without it being a true mathematical relation, i.e., roughly $10^{-1000}$. We therefore believe that our conjectures are mathematical formulas awaiting formal proof. Importantly, this estimate of likelihood of course does not substitute the need for formal proof. We hope that the discovered conjectures could act as an accelerator for the discovery of new mathematical structures and aid in future discoveries.

Apart from the algorithm, the analysis and manipulation of CFs in our work relies on several novel results that are presented below. We present a connection between simple CFs with interlaced polynomial sequences and the SICF structure found in our results. Furthermore, we present the **Folding Transform**, which reveals a connection between polynomial CFs and CFs made of interlaced polynomial sequences or more generally polynomial matrices. These connections provide insight on the space of constants for which we can expect the ESMA algorithm to find formulas.

# 2 Background

## 2.1 The Berlekamp-Massey algorithm in ESMA

We present a new application of the Berlekamp-Massey algorithm [8-9], used for pattern recognition in an automated search of conjectures on mathematical constants. Given an integer sequence, the Berlekamp-Massey algorithm finds the minimal linear recurrence with integer coefficients that can produce the sequence, returning the coefficients of the recurrence relation [9].

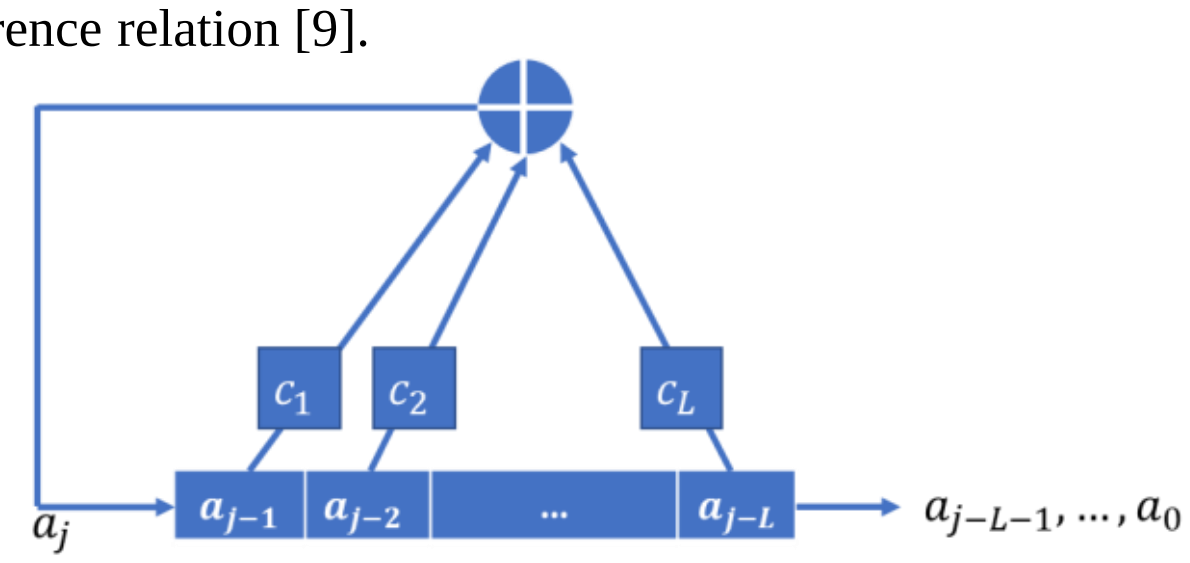

**Figure 1 | Linear-Feedback Shift Register (LFSR)** containing $L$-cells, where $a_{j-k}$ is the value held in the $k^{th}$ cell, and $c_k$ is the coefficient with which we multiply the value before linear feedback.

The Berlekamp-Massey algorithm is based on Berlekamp's decoding algorithm [9,12], generalized by James L. Massey to solve the shift register synthesis problem, the task of finding the shortest **Linear-Feedback Shift Register (LFSR)** that outputs a given sequence over a finite field [12-13]. In such registers, the input at every clock is a fixed linear function of the register's current state, creating a linear recursion with coefficients referred to as connection coefficients [9] (see Figure 1). Finding the shortest integer coefficient recurrence relation of a sequence is analogous to finding the shortest linear-feedback shift register that can output the sequence, and this is what the Berlekamp-Massey algorithm achieves.

A linear-feedback shift register of length $L$ with the initial contents of the $L$ cells given by $a_0, a_1, \dots a_{L-1}$ (initial conditions) and connection coefficients given by $\{c_i \in \mathrm{GL}(p) | i \in \{1,2,\dots,L\}\}$, for some prime number $p$, produces the following output sequence over $\mathrm{GL}(p)$ [9]:

$$\left(a_j + \sum_{i=1}^{L} c_i a_{j-i}\right) \bmod \mathrm{p} = 0 \quad \forall j \in \{L, L+1, \dots, n-1\} \tag{5}$$

Given a sequence $a_0, \dots, a_n$, the Berlekamp-Massey algorithm finds and returns the connection coefficients $c_i$ of the minimal linear-feedback shift register for which (5) is satisfied; these are coefficients of the recursive relation that describes the sequence.

The length of the resulting recursion, which is the number of cells, further reveals to us the significance of the pattern detected. Given an input sequence of length $n$, if the resultant recurrence is of length $L \geq n/2$, then it is a trivial solution. In this case, we can simply take the initial conditions of the LFSR to be the first elements of the sequence and calculate the connection coefficients so they create the next $n/2$ elements. We therefore add a verification step by using the length of the register found by the algorithm to identify whether the pattern found in the sequence is significant and unique, requiring a recurrence length $L < n/2$ (more details in [12]). Even then, there is still a probability for a false positive that scales as $\sim(1/p)^{n-2L}$, and therefore we apply a second verification stage.

In ESMA, we apply the Berlekamp-Massey algorithm on $a_j$ sequences extracted from the expansion of CFs as shown in the next sections. We use a large finite field $\mathrm{GL}(p)$ so that the connection parameters are found explicitly (rather than modulo the finite field). Specifically, the presented results were found with $p = 199$. One can use an even smaller value of $p$ to make the algorithm run faster, at the price of post-processing for extracting the explicit connection parameters from their values in the finite field. For conjecture verification, we use the extracted recurrence relation to generate the sequence over $\mathbb{Z}$ to any desired length (not limiting the generated numbers to GL(p)). Altogether, ESMA enables us to identify any pattern that can be described as an integer coefficient recurrence relation.

## 2.2 Signed Interlaced Continued Fractions

In the most general sense, **Signed Interlaced Continued fractions (SICFs)** are CFs in which the partial denominator is an interlaced sequence, and the partial numerator is some periodic sequence of $\pm 1$, $b_j \in \{-1,1\}$, whose period, $\beta$, denotes the number of sub-sequences that are interlaced. For example, the following conjecture found on Bessel functions of the first kind, $\mathrm{J}_y(x)$ (of order $y$ at point $x$), by the ESMA algorithm is an example of an SICF with $a_j$ made of interlaced linear sequences:

$$\frac{2}{\frac{\mathrm{J}_5(1)}{\mathrm{J}_3(1)}+1} = 2 - \frac{1}{40 - \frac{1}{3 - \frac{1}{56 - \frac{1}{4 - \cdots}}}},\ b_j = -1,\ a_j = \begin{cases} 2+k & j = 2k \\ 40+16k & j = 2k+1 \end{cases}. \tag{6}$$

$$k \in \{0,1,2,\ldots\}, \qquad \forall j \in \mathbb{N}$$

SICFs are the structure with which ESMA extracts conjectures, where the interlaced sequences are given by any sequence that can be modeled by a linear-feedback shift register. In the above example the $a_j$ sequence is described by the following integer coefficient recurrence relation:

$$a_j - 2a_{j-2} + a_{j-4} = 0,$$

with initial conditions:

$$a_0 = 2, a_1 = 40, a_2 = 3, a_3 = 56.$$

In all the results so far, we obtain an $a_j$ made of positive interlaced linear sequences ($a_j > 0\ \forall j \in \mathbb{N}$). Therefore, for the rest of the paper when referring to SICF, we refer specifically to those with positive polynomial sequences in the partial denominator, satisfying $a_j > 0\ \forall j \in \mathbb{N}$. To enable discovery of conjectures in the form of an SICF, we introduce a method to extract sample of the $a_j$ sequence for CFs with a signed partial numerator, $b_j = \pm 1$.

## 2.3 Extraction of Signed Interlaced Continued Fractions and the Euclidean algorithm

We generalize the conventional method of calculating the CF of a constant [14-16, 22] to enable sign variation in the partial numerator. This procedure enables us to expand every constant to a SICF with any sequence of $\pm$ signs in the partial numerators, allowing for a larger space of candidate formulas to be analyzed.

Simple CFs are unique, every irrational number has a single simple CF which it is equal to, and every rational number has 2 simple CF expansions [14]. Calculating the simple CF of some constant $\alpha$ is simply an application of a non-terminating Euclidean Algorithm with $\alpha$ and 1 [14], where the respective quotients form the $a_j$ sequence and $b_j = 1\ \forall\ j \in \mathbb{N}$. When enabling sign variation in the Euclidean algorithm, we provide a degree of freedom. Sign variation offers more options in each iteration of the Euclidean algorithm by allowing for negative remainders. In many cases, sign variation enables the algorithm to terminate in less iterations, meaning a more efficient representation is found. In terms of CFs, sign variation ($b_j = \pm 1\ \forall j \in \mathbb{N}$) allows for the extraction of additional CF representations for a single constant, some of which may have more favorable properties. Take a simple example of finding the $\gcd(14, 9)$ and the resultant CF formed:

| | **Without sign variation** | **With sign variation** | |
|---|---|---|---|
| | $14 = 9 \times \mathbf{1} + 5$ | $14 = 9 \times \mathbf{2} - 4$ | |
| $\frac{14}{9} = 1 + \frac{1}{1+\frac{1}{1+\frac{1}{4}}}$ | $9 = 5 \times \mathbf{1} + 4$ | $9 = -4 \times -\mathbf{2} + 1$ | $\frac{14}{9} = 2 - \frac{1}{2+\frac{1}{4}}.$ |
| | $5 = 4 \times \mathbf{1} + 1$ | $-4 = -1 \times \mathbf{4} + 0$ | |
| | $1 = 1 \times \mathbf{1} + 0$ | | |

We see that enabling signed variation results in a more efficient calculation, 3 iterations rather than 4, along with a reduced CF expansion of the number. The signed CF provides us with a larger number of possibilities to search for patterns. As a result, we enable the discovery of more conjectures and occasionally discover more efficient expansions. For example, consider the conjecture on $J_5(1)/J_3(1)$, seen in (6):

$$\frac{2}{\frac{J_5(1)}{J_3(1)}+1} = 2 - \cfrac{1}{40-\cfrac{1}{3-\cfrac{1}{56-\cfrac{1}{4-\cdots}}}},\ b_j = -1,\ a_j = \begin{cases} 2+k & j = 2k \\ 40+16k & j = 2k+1 \end{cases} \quad \forall j \in \mathbb{N}, k \in \{0,1,2,3 \dots\}$$

The new and unproven conjecture converges to $J_5(1)/J_3(1)$ at a rate of 4.3008 digits per term (averaged over 100 terms), whilst the simple CF expansion of $J_5(1)/J_3(1)$ converges at a rate of 1.9046 digits per term (see Figure 4 in Appendix A).

In many cases, the SICF that results from ESMA is simpler than the simple CF of the same constant. Take for example the simple CF of $(2+2\mathrm{e})/(-1+3\mathrm{e})$. It's simple CF involves an interlaced sequence of period 8 (up to a couple of initial elements). In comparison, the SICF is of period 6. Equation (6) presents a similar example for $J_5(1)/J_3(1)$, whose simple CF is of period 16 whereas its SICF is of period 2 (up to initial elements). Therefore, the SICF provides us with a larger search space where we can detect simpler patterns and, in some cases, even faster converging expansions (Figure 4) than with the simple CF counterpart.

We additionally find an algorithmic advantage in utilizing SICF as more conjectures for the same search space are found. For the same search space of coefficients, coefficients between -3 and 3 and polynomial degree of 1, there was an increase of ~357% in the number of conjectures found on e when searching with SICFs rather than only simple CFs. A similar increase in the number of conjectures (as a result of enabling sign variation) is observed across all other constants. There was an increase of ~591% in the number of conjectures on $\tan(1)$ in same above search space. Presenting the algorithmic advantage gained by signed CF extraction.

<table>
<tr>
<td>$a_0 = \lfloor c \rfloor\ if\ (b_1 > 0){:}\ else\ \lceil c \rceil$<br>$c = c - a_0$<br>$for\ i = 1{:}\ depth{:}$<br>$if c = 0\ return$<br>$c = \frac{b_i}{c}$<br>$a_0 = \lfloor c \rfloor\ if\ (b_{i+1} > 0){:}\ else\ \lceil c \rceil$<br>$c = c - a_i$<br>Signed CF Extraction algorithm</td>
<td rowspan="2">$b = [b_1, b_2, b_3, \dots] = [-1,1,1,\dots]$<br>$c = \frac{2+2\mathrm{e}}{-1+3\mathrm{e}} = 1.039374$<br>$b_1 = -1 < 0 \ \rightarrow a_0 = \lceil c \rceil = \mathbf{2}$<br>$c = c - a_0 = -0.960626$<br>$c = \frac{b_1}{c} = \frac{-1}{c} = 1.04098$<br>$b_2 = 1 > 0 \ \rightarrow a_1 = \lfloor c \rfloor = \mathbf{1}$<br>$c = c - a_1 = 24.402$<br>$c = \frac{b_2}{c} = \frac{1}{c} = 1.04098$<br>$b_3 = 1 > 0 \ \rightarrow a_2 = \lfloor c \rfloor = \mathbf{24}$<br>...<br>Example Extraction of Signed CF formula</td>
</tr>
<tr>
<td>$\frac{2+2\mathrm{e}}{-1+3\mathrm{e}} = \mathbf{2} - \cfrac{1}{\mathbf{1} + \cfrac{1}{\mathbf{24} + \cfrac{1}{\dots}}}$<br>Example Signed CF formula</td>
</tr>
</table>

To extract a signed CF, we accommodate for having some $b_j = -1$ and we extract the integer part of the CF using ceil operator rather than the floor operator used for $b_j = 1$ (as in the Euclidean Algorithm). For example, for some constant $c$ the extraction of its SICF can be seen in the table above. We utilize each given signed $b_j$ sequence to extract a sample of the $a_j$ integer sequence directly from the decimal representation of a constant. Consequently, we obtain an integer sequence for which we can attempt to recognize a pattern using the Berlekamp-Massey algorithm.

# 3 The Extract Signed-CF Massey-Approve (ESMA) Algorithm

We present a novel algorithm (Figure 2) that extracts a Signed Interlaced Continued Fraction (SICF) of a given constant $c$ in the following form:

$$\frac{f_{m,L}(c)}{g_{m,L}(c)} = a_0 + \mathbb{K}_1^\infty \frac{b_j}{a_j} = a_0 + \cfrac{b_1}{a_1 + \cfrac{b_2}{a_2 + \cfrac{b_3}{a_3 + \cdots}}}. \tag{7}$$

Where $f_{m,L}$ and $g_{m,L}$ are integer polynomials whose degrees are at most $m$ with coefficients over a range of integer values $[-L, L]$. In this case, $a_0 + \mathbb{K}_1^\infty b_j/a_j$ is some SICF with partial numerator $b_j \in \{1, -1\}$ of maximal period $\beta_b$ and partial denominator $a_j \in \mathbb{Z}$.

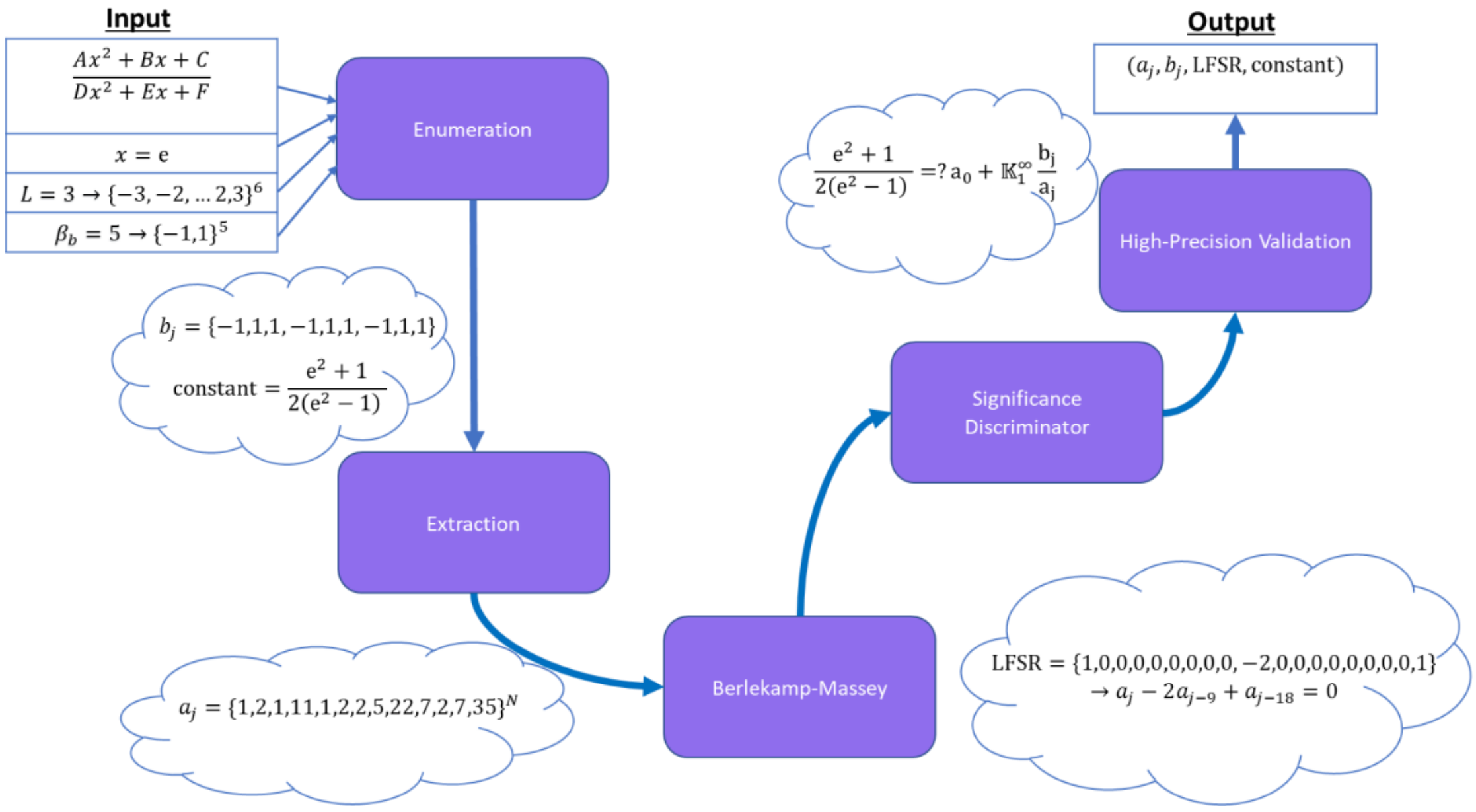

**Figure 2 | The Extract Signed-CF Massey Approve (ESMA) Algorithm:** given a fundamental constant, a polynomial degree, range of coefficient values $L$, and a maximal period for signed $b_j$ sequences $\beta_b$, we enumerate over all possible rational functions and over all possible $b_j$ sequences. For each non-trivial expression, we extract its CF with each $b_j$ sequence to a low depth and then utilize the Berlekamp-Massey algorithm to detect pattens in extracted sample of the $a_j$ sequence. Results are finally validated to a high precision (1000 decimal places) and then stored as conjectures. The "clouds" in the figure present an example run for a single constant value and $b_j$ pair.

We begin by enumerating over all possible non-trivial rational functions $f_{m,L}(x)/g_{m,L}(x)$ (see Figure 2 for an example). For each rational function, we substitute the mathematical constant $c$ (evaluated to 1000 decimal places in the examples shown in this work). We then enumerate over all periodic $b_j$ sequences with periods between 1 and $\beta_b$ formed of $\pm 1$. For each $b_j$ sequence and constant value pair, an $a_j$ sequence is extracted up to a finite depth $N$, using the extraction algorithm introduced earlier. The sequence is passed to the Berlekamp-Massey algorithm in attempt to recognize a significant recurrence pattern, i.e., a minimal length linear feedback-shift register. If the resulting register length is shorter than half the length of the extracted sequence (see section 2.1), it is considered significant and saved for verification. The resulting pair of $a_j$, $b_j$ sequences represent a SICF.

We verify the result by utilizing the obtained recursion parameters to calculate $a_j$ to a greater depth and then evaluate the SICF to compare it with the constant. The SICF is efficiently calculated by utilizing the recursive formula for numerators $p_j$ and denominators $q_j$ [17]:

$$p_j = a_j p_{j-1} + b_j p_{j-2}\ ,\ q_j = a_j q_{j-1} + b_j q_{j-2} \tag{8}$$

$$p_{-1} = 1, \qquad p_0 = a_0, \qquad q_{-1} = 0, \qquad q_0 = 1,$$

yielding $p_j/q_j$ as a rational approximation of the constant given by the first $j$ elements of the partial numerator and denominator sequences of the SICF. Each case for which the numerical values are identical for up to 1000 decimal places is considered a new conjecture. For the full implementation of the algorithm, refer to our git https://github.com/RamanujanMachine/RamanujanMachine/tree/master/ESMA.

Later additions to the algorithm enable users to construct custom made function generators easily, not limiting the search to rational polynomial functions $f_{m,L}(c)/g_{m,L}(c)$ but rather to any family of parametric functions over a discrete parameter space.

The computational complexity of the algorithm depends on the space over which we enumerate. The space of rational functions $f_{m,L}(x)/g_{m,L}(x)$ is $O\left((2L+1)^{2(m+1)}\right)$. For convenience, our system supports saving symbolic enumerations locally so we can substitute different constants into the symbolic expressions. The rational function enumerations are also simplified to reduce redundancies, ensuring no trivial cases are stored (rational numbers, where polynomials cancel out).

The $b_j$ sequence enumeration takes $O(2^{\beta_b})$, where $\beta_b$ is the maximal period. In our most common searches, $\beta_b$ ranges between 1 and 5, and thus we enumerated $2^1 + 2^2 + 2^3 + 2^4 + 2^5$ combinations of $\pm 1$ sequences. For each $b_j$ sequence and constant value pair, we extract the $a_j$ sequence to a depth $N$ (the results presented in this work are found using $N = 50$). The complexity of this extraction is negligible relative to the application of the application of the Berlekamp-Massey algorithm, whose complexity is $O(N^2)$ [12, 18]. In case that a result was found, we verify it to 1000 decimal places by evaluating the resulting SICF and comparing with the original constant.

Given the above considerations, the time complexity of the algorithm is given by $O\left((2L+1)^{2(m+1)} 2^{\beta_b} N^2\right)$. While this may seem computationally expensive, it has two substantial advances over exhaustive search method such as the first Ramanujan Machine algorithm (meet-in-the-middle regular formulas) in [6]: (1) The enumeration over possible polynomial CFs and the expansion of each one, which are the most expensive operations, are replaced by enumerating over the signed sequences $O(2^{\beta_b} N^2)$. (2) Even the enumeration over the rational functions is more efficient than in [6], since infinite many cases of rational functions are captured by a single run of the Berlekamp-Massey algorithm, inside the initial conditions in the register. For example, both the expressions $x$, $1/x$, $1/x + k$, $x + k$, and many other Mobius transforms of $x$, are covered by the same instance of running the Berlekamp-Massey algorithm, as initial sequence elements which do not follow the found pattern can be simply absorbed to the initial conditions of the recurrence relation. This applies to any value of $x$ including the rational functions of $x$, presenting a more efficient enumeration.

However, despite its efficiency, the disadvantage of ESMA relative to the first Ramanujan Machine algorithm [6] lies in the space of CFs that can be discovered. This space is analyzed in the following sections.

# 4 The Space of Constants Captured by the ESMA Algorithm

## 4.1 Example Results

The proposed ESMA algorithm has discovered many previously unknown conjectures, some having faster converging expansions for various constants (see Appendix A, Figure 4) relative to their known simple CFs. In Figure 3, we see a sample of conjectures and known formulas found by ESMA, all of which converge at a super-exponential rate (error decreasing as $\approx e^n/n!$ for CF depth $n$ [6]). This rate of convergence is shown in Figure 3, compared with the (slower) exponential convergence of the Golden Ratio's simple CF. Similar fast rates of convergence are seen in most conjectures found by ESMA. To better understand our results and the reason for their super-exponential convergence rate in most cases, we analyze the representable set of the SICF structure over which ESMA searches.
Additional example results are presented in Appendix A.

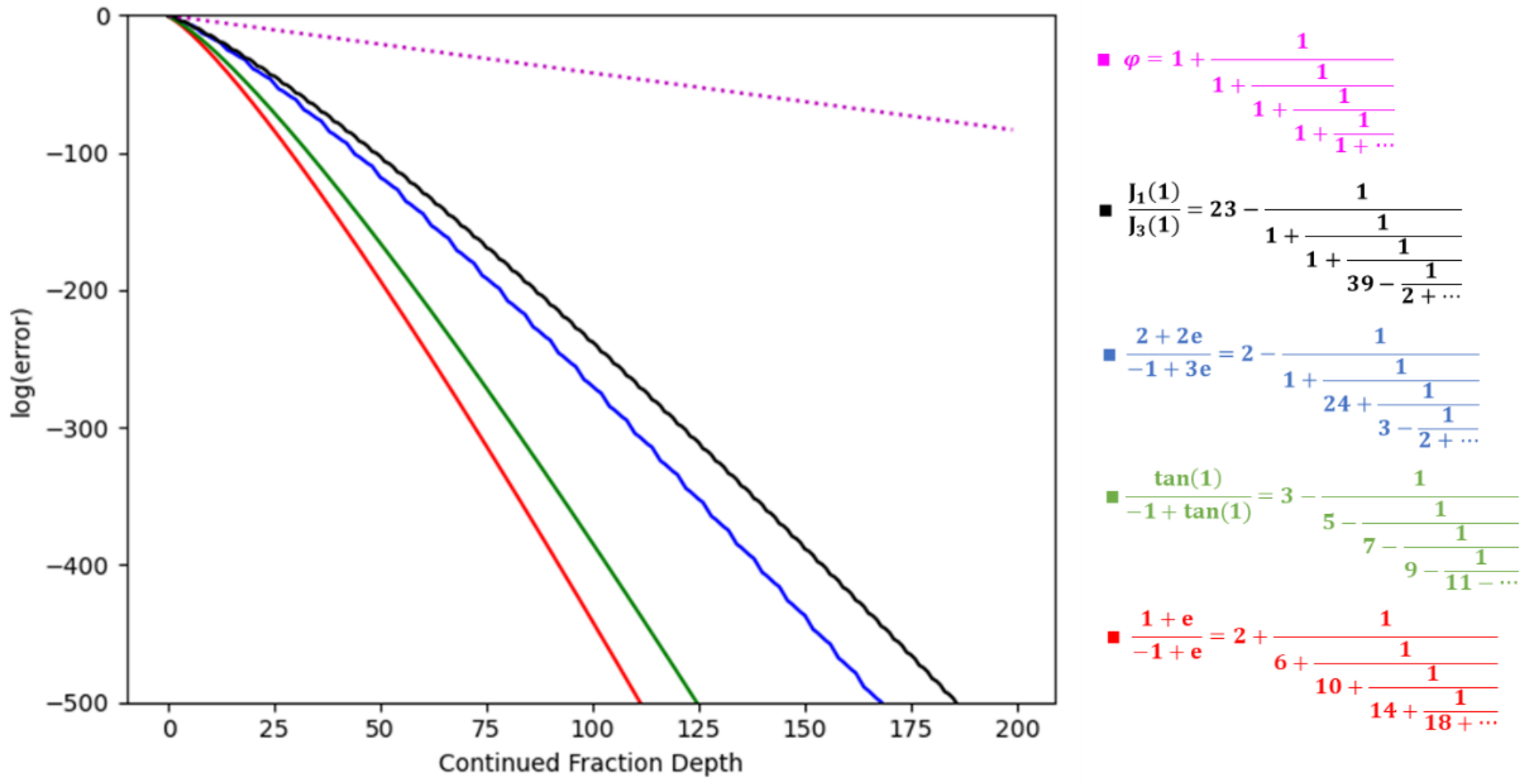

**Figure 3 | Convergence rates of the SICF conjectures**. The figure presents the approximation error, the log of the absolute difference between the SICFs approximation at a given CF depth and the fundamental constant ($log_{10}(error)$ vs SICF depth). All the above results converge super-exponentially (validated numerically), the exponentially converging formula for $\varphi$ is plotted for reference to visualize the difference between exponential and super-exponential convergence. The red, green, and pink lines describe the convergence of known formulas, the others of new and unproven conjectures.

## 4.2 Equivalent Representations of Constants

As described in the previous section, the ESMA algorithm searches for a CF expansion of a given constant where the partial numerator $b_j$ is periodic with $\pm 1$ entries and the partial denominators $a_j$ satisfy some recurrence relation, $a_0 + \mathbb{K}_1^\infty b_j/a_j$. While a linear recurrence determines the sequence completely given "enough" initial conditions, we can sometimes decompose it into several subsequences and apply the recurrence to each one of them separately. For example, with the linear recurrence $a_j - 2a_{j-2} + a_{j-4} = 0$ mentioned earlier in section 2, we can apply the recurrence to the subsequence at even and odd indices separately obtaining an interlaced sequence. Thus, we can think of the sequence as $a_j = \begin{cases} A_1(n) & j = 2n-1 \\ A_2(n) & j = 2n \end{cases}$ (if we ignore the $a_0$ element). We now have a natural decomposition of the $b_j$ and $a_j$ sequences to $\beta_b$ and $\beta_a$ subsequences respectively, and up to taking the least common multiplier of $\beta_b$ and $\beta_a$ we may assume that both are equal to the same $\beta$, which we call the **period** of the SICF. Hence, we have functions $B_i,\ A_i: \mathbb{N} \to \mathbb{Z}$ for $i = 1, \dots, \beta$ such that $b_{(n-1)\times\beta+i} = B_i(n)$, $a_{(n-1)\times\beta+i} = A_i(n),\ \forall n \in \mathbb{N}\backslash 0$.

While $B_i$ are always constant 1 or $-1$, in general, ESMA may find many types of $A_i$ functions (though they must be some combination of polynomials and exponents). However, in all results so far, the functions found were linear polynomials (see Appendix A), leading us to question what sort of numbers can be represented with polynomial $A_i$'s and in a more general sense what constants one can expect the ESMA algorithm to catch (find conjectures for). We denote this set of numbers, the representable set of the SICF structure, as:

$$\mathfrak{R}_1 = \left\{ \mathbb{K}_1^\infty \frac{b_j}{a_j} \;\middle|\; \begin{array}{c} \exists \beta \text{ and } A_i, B_i\ \ i = 1, \dots, \beta, \text{ s.t.} \\ b_{(n-1)\times\beta+i} = B_i(n)\ \ i \in \{1, \dots, \beta\}\ \ B_i \equiv \pm 1 \\ a_{(n-1)\times\beta+i} = A_i(n)\ \ i \in \{1, \dots, \beta\}\ \forall n \in \mathbb{N}\backslash 0\ \ A_i \text{ is polynomial} \\ \mathbb{K}_1^\infty \frac{b_j}{a_j} \text{ converges} \end{array} \right\}. \tag{9}$$

The above set refers to the set of all possible numbers that have a formula with SICF up to some Mobius transform. Similar to common standard CF representations [23], here too we can represent the more general CF using Mobius transforms (for details, see Appendix B.1). More specifically, given a Mobius transform $\begin{pmatrix} a & b \\ c & d \end{pmatrix}(z) = \frac{az+b}{cz+d}$ $a, b, c, d \in \mathbb{Z}$, we get that

$$\mathbb{K}_1^\infty \frac{b_j}{a_j} = \lim_{N\to\infty} \prod_{j=1}^N \begin{pmatrix} 0 & b_j \\ 1 & a_j \end{pmatrix}(0) = \lim_{N\to\infty} \prod_{j=1}^{N+1} \begin{pmatrix} 0 & b_j \\ 1 & a_j \end{pmatrix}(\infty). \tag{10}$$

In our interlaced presentation above, the matrices come in natural batches of size $\beta$, and it is only natural to multiply each such batches together to form a single polynomial matrix, which we refer to as the **collapsed matrix**:

$$\prod_{j=1}^{N\beta} \begin{pmatrix} 0 & b_j \\ 1 & a_j \end{pmatrix} = \prod_{n=1}^{N} \left[ \prod_{i=1}^{\beta} \begin{pmatrix} 0 & b_{(n-1)\times\beta+i} \\ 1 & a_{(n-1)\times\beta+i} \end{pmatrix} \right],\ M_n = \prod_{i=1}^{\beta} \begin{pmatrix} 0 & b_{(n-1)\times\beta+i} \\ 1 & a_{(n-1)\times\beta+i} \end{pmatrix}, n \in \mathbb{N}\backslash 0. \tag{11}$$

We denote the collapsed matrix by $M_n$ and we get that the entries of $M_n$ are polynomial in $n$. The collapsed matrix can additionally be used to represent the more general case where both our $a_j$ and $b_j$ sequences are interlaced polynomial sequences, an **interlaced CF**. Thus, we can automatically deduce that $\mathfrak{R}_1$ above is contained in a much more general set,

$$\mathfrak{R}_2 = \left\{ \lim_{N\to\infty} \left[\prod_{n=1}^N M_n\right](\infty) \;\middle|\; M_n \text{ is a } 2\times 2 \text{ matrix with polynomial entries} \right\}. \tag{12}$$

While the presentation as a product of a polynomial matrix can seem very general, we can show that it can almost always be presented in a much simpler manner, namely as a polynomial CF. Polynomial CFs are well known [6] [7, 19, 20] and have been shown to enable swift proofs on the properties of CFs and the constants they can represent, thus providing information on a seemingly more general structure (e.g., Apery's proof of $\zeta(3)$'s irrationality).

**Theorem 1.** Let $M_n = \begin{pmatrix} c_n & d_n \\ e_n & f_n \end{pmatrix}$, satisfying $e_n \neq 0$, be some polynomial matrix ($c, d, e, f \in \mathbb{Z}[x]$) for which the following limit exists $\lim_{N\to\infty} \left[\prod_{n=1}^N M_n\right](\infty)$. Then there exists some Mobius transform $T$ with entries in $\mathbb{Z}$, and polynomials $a', b' \in \mathbb{Z}[x]$ such that:

$$\lim_{N\to\infty}\left[T\prod_{n=1}^{N}M_n\right](\infty) = \mathbb{K}_1^{\infty}\frac{b'(n)}{a'(n)}$$

Proof. See Appendix B.1

Note that since the Mobius transform $T$ has its entries in $\mathbb{Z}$, a number $\alpha$ is rational if and only if $T(\alpha)$ is rational. Hence, up to Mobius equivalence we see that $\mathfrak{R}_2$ is contained in the set of numbers with expansions as polynomial CFs,

$$\mathfrak{R}_3 = \left\{\mathbb{K}_1^{\infty}\frac{b'(n)}{a'(n)}\middle| a', b' \in \mathbb{Z}[x]\right\}. \tag{13}$$

As a result, we can deduce that up to an integer Mobius transformation $\mathfrak{R}_1 \subseteq \mathfrak{R}_2 \subseteq \mathfrak{R}_3$.

We present the **Folding Transform** $T$, converting any polynomial matrix satisfying the theorem's conditions to a polynomial CF (see Appendix B.1). In other words, the Folding transform can convert any constant equal to an interlaced CF to a polynomial CF, up to some Mobius transform. For any $\alpha \in \mathbb{R}$ equal to a general interlaced CF where both $b_j$ and $a_j$ are interlaced polynomial sequences, and specifically for a SICF where $b_j = \pm 1$, we can multiply all matrices in each period to obtain a polynomial matrix, $M_n$. Applying the Folding transform on the resulting polynomial matrices, we obtain a polynomial CF that we denote by $T(\mathrm{M_n})$. We say that $\alpha$ and $T(M_n)$ are **semi-equivalent**. For example, applying the Folding transform on the simple CF of e, we obtain the following polynomial CF (for the full derivation see Appendix B.3):

$$\mathrm{e} - 2 = \frac{1}{1+\frac{1}{2+\frac{1}{1+\frac{1}{1+\frac{1}{4+\cdots}}}}} \rightarrow^{T} \frac{-8\mathrm{e}+19}{\mathrm{e}-2} = \frac{4\times 2^2+4\times 2-3}{8\times 2^2+16\times 2+8+\frac{4\times 3^2+4\times 3-3}{8\times 3^2+16\times 3+8+\frac{4\times 4^2+4\times 4-3}{8\times 4^2+16\times 4+8+\cdots}}} \tag{14}$$

Note that the proof for the formula on the right can be directly deduced from the Folding transform process and using e′s standard CF. The Folding transform allows us to better understand our largely unexplored SICF structure by utilizing prior knowledge of polynomial CF properties to understand what kind of constants we can represent. To further narrow down our analysis, we present a novel connection showing that constants with SICF expansions can be presented in simpler forms.

## 4.3 Representable set of Signed Interlaced Continued Fractions (SICFs)

At first glance, the representable set of SICFs, $\mathfrak{R}_1$, may seem to expand the set of numbers we can represent with a simple CF by allowing for sign variation. However, we find that every SICF expansion of an irrational constant can be presented as a simple CF with a "regular" pattern up to some Mobius mapping. Rational numbers may have SICF expansions but always have a finite CF and are not of interest to us. We denote the representable set of simple CFs with interlaced polynomial sequences, also referred to as **simple interlaced CFs**, as:

$$\mathfrak{R}_4 = \left\{\mathbb{K}_1^{\infty}\frac{1}{a_j}\middle| \begin{array}{c} \exists\beta \text{ and } A_i, B_i \;\; i = 1, \dots, \beta, \;\; s.t. \\ a_{(n-1)\times\beta+i} = A_i(n) \;\; i \in \{1, \dots, \beta\} \;\; \forall n \in \mathbb{N}\backslash 0 \;\; A_i \text{ polynomial} \end{array}\right\}. \tag{15}$$

Unlike SICFs, simple CFs have been studied extensively [7, 14, 17, 21, 22], much is known about their properties and about constants' simple CF expansions. Therefore, discovering that for irrational numbers $\mathfrak{R}_1 \equiv \mathfrak{R}_4$ up to some Mobius transform shines light on what constants we can expect ESMA to catch: those whose simple CF expansions have a pattern in their partial denominator. In other words, we expect ESMA to catch numbers who have simple interlaced CF expansions.

**Theorem 2.** For any signed interlaced continued fraction equal to $\alpha \in \mathbb{R}/\mathbb{Q}$, there exists some Mobius transform of $\alpha$ which is equal to a simple interlaced continued fraction.

Proof. See Appendix C.1

The proof of Theorem 2 is constructive, meaning we have an algorithm to convert every SICF to a simple interlaced CF expansion. For example, taking the following conjecture on $\frac{2}{\tan(1)}$:

$$\frac{2}{\tan(1)} - 2 = -\frac{1}{2-\frac{1}{2-\frac{1}{3-\frac{1}{14-\cdots}}}},\;\; b_j = -1,\;\; a_j = \begin{cases} 3n-1 & j = 4n-3 \\ 2 & j = 4n-2 \\ 3n & j = 4n-1 \\ 2+12n & j = 4n \end{cases}, \tag{16}$$

$$\forall j, n \in \mathbb{N}\backslash 0.$$

We use the method shown in the proof of Theorem 2 to obtain its simple interlaced CF form up to a Mobius transform (for full derivation see Appendix C.2):

$$-\frac{1}{\frac{2}{\tan(1)}-2}-1=\frac{1}{2+\frac{1}{1+\frac{1}{12+\frac{1}{1+\frac{1}{3+\cdots}}}}}\ ,\ b_j=1,\ a_j=\begin{cases}3n-2 & j=6n-5\\ 1 & j=6n-4\\ 12n & j=6n-3\\ 1 & j=6n-2\\ 3n & j=6n-1\\ 2 & j=6n\end{cases}, \tag{17}$$

$$\forall j,n\in\mathbb{N}\backslash 0.$$

As a result of Theorem 2, we can fully characterize the representable set of SICFs, $\mathfrak{R}_1$, with simple interlaced CFs, a much less general structure. Through this connection we reveal a limiting property of the ESMA algorithm, only constants with regular patterns or interlaced sequences in the partial denominator of their simple CF (which is unique) can be caught by ESMA. This set of numbers includes e, $\tan(1)$, second-degree algebraic numbers (see Appendix B), and more. The connection found further explains why constants such as π were not be caught by the algorithm, as there is no discernable pattern in π's simple CF [24]. We utilize this theorem to simplify our analysis of SICFs. We look at polynomial CF expansions semi-equivalent to simple interlaced CFs, who have less degrees of freedom than SICFs as $b_j=1\ \forall j\in\mathbb{N}\backslash 0$, and analyze their properties to better characterize the representable set of SICFs.

We find that the representable set of simple interlaced CFs with non-constant partial denominator sequences is characterized by irrational constants whose polynomial CF expansions converge super-exponentially revealing to us why most of our results are of constants who have super-exponentially converging expansions. We do not analyze simple interlaced CFs with constant partial denominator sequences as they are trivially characterized by the 2$^{nd}$ algebraic numbers [25-26].

**Theorem 3.** Given a simple interlaced continued fraction satisfying $B_i(n)=1\ \&\ A_i(n)>0$, where $A_i,B_i\in\mathbb{Z}[\mathrm{x}]$, $i\in\{1,..,\beta\},\forall n\in\mathbb{N}$ of period $\beta$, where $\exists i\in\{1,..,\beta\}$ s.t $\deg(A_i)>0$, its semi-equivalent polynomial continued fraction's partial numerator $b'_n$ and partial denominator $a'_n$ ($\forall n\in\mathbb{N}$) have the following polynomial degrees:

$$\deg(b')=\sum_{i=1}^{\beta-1}2\times\deg(A_i),\qquad \deg(a')=\left[\sum_{i=1}^{\beta-1}2\times\deg(A_i)\right]+\deg(A_\beta)$$

Note, $A_i$ is the $i^{\text{th}}$ polynomial sequence of $a'_j$s $\beta$ different sub-sequences

Proof. See Appendix D.1

Through the utilization of known properties of polynomial CFs [7, 19-20] the above theorem provides us with valuable information to better understand the rate of convergence of these polynomial CF representations and the irrationality of the constants ESMA can catch.

**Corollary 2.** A polynomial continued fraction semi-equivalent to a simple interlaced continued fraction, satisfying $\exists i\in\{1,..,\beta\}$ s.t $\deg(A_i)>0$, converges super-exponentially.

Proof. See Appendix D.2

**Lemma 3.** The partial denominator sequence of a polynomial continued fraction semi-equivalent to a simple interlaced continued fraction is positive: $a'_n\in\mathbb{Z}[\mathrm{x}],\ a'_n>0\ \forall n\in\mathbb{N}\backslash 0$

Proof. See Appendix D.1

**Corollary 3.** Any number $\alpha\in\mathbb{R}$ equal to a simple interlaced continued fraction converges to an irrational limit.

Proof. See Appendix D.3

As can be directly deduced from Theorem 3 and Lemma 1, the set of numbers representable by a simple interlaced CF, $\mathfrak{R}_4$, is a subset of irrational numbers that have polynomial CF expansion which converge super-exponentially. Specifically, those with a positive partial denominator whose degree is greater than or equal to that of the partial numerator, $\mathfrak{R}_4\subset\{\mathfrak{R}_3|\deg(b')\le\deg(a'),a'(n)>0\forall n\in\mathbb{N}\backslash 0\}$. Recall, $\mathfrak{R}_3$ refers to the representable set of polynomial CFs. We can therefore deduce, that as $\mathfrak{R}_1\equiv\mathfrak{R}_4\subset\mathfrak{R}_3$, meaning numbers who have SICF expansions are a subset of numbers which have polynomial CF expansions.

Overall, the representable set of SICFs, which is also the space of constants we expect ESMA to catch, is characterized by irrational constants whose simple CF expansion is made of interlaced polynomial sequences (up to a Mobius transform). Such CFs with constant partial denominator sequences converge at an exponential rate to a $2^{nd}$ degree algebraic numbers. Such CFs with non-constant sequences have polynomial CF expansions that converge at a super-exponential rate. Intriguingly, ESMA could potentially find CFs with non-polynomial sub-sequences. However, our extensive runs on ESMA with up to 50 different constants only found polynomial sub-sequences and thus further supports the conclusions of our mathematical analysis.

Prior to our mathematical analysis, we only obtained results on e. Knowing the properties of the polynomial CF forms of constants ESMA can catch, we specifically targeted and searched the literature for constants with those properties. As a result of this targeted search, we found various novel conjectures on constants such as the ratio of Bessel functions at different points, $\tan(1)$, $\tanh\left(\frac{1}{4}\right)$, and more (Appendix A).

# 5 Discussion and Outlook

The ESMA algorithm proved successful in discovering various conjectures on constants that converge significantly faster than their simple CF expansions. For example, the algorithm discovered conjectures on $J_1(1)/J_3(1)$, $J_5(1)/J_3(1)$, and the Golden Ratio $\varphi$, which converge faster than their simple CFs. A conjecture found on $\varphi$ converges approximately 6 times faster than its simple CF and on $J_5(1)/J_3(1)$ more than twice as fast (see Appendix A Figure 4). In Appendix A, we plot the convergence of a sample of faster converging conjectures alongside their simple CF expansions, presenting the significant improvement in the convergence rate of such expansions.

Furthermore, ESMA comes as a complementary approach to the existing exhaustive search methods for automation of mathematical discovery, such as the Ramanujan Machine algorithm in [6]. ESMA improves the efficiency of enumeration over the space of possible CFs and rational functions.

Despite its advantages, ESMA also has important shortcomings. ESMA seems unable to find conjecture formulas for certain constants. Specifically, all formulas found by ESMA (apart from trivial ones) have a super-exponential convergence rate, and it is thus unable to find formulas that converge more slowly. This limitation implies that ESMA cannot find conjectures for constants for which no super-exponentially converging formula exists. The limitation is even more severe – any constant which ESMA can catch must also have an interlaced polynomial pattern in its simple CF. Such a pattern is not known for mathematical constants such as $\pi$, $\zeta(3)$, and $G$, and thus we expect ESMA cannot find conjectures for them, unless such a formula could be found for a transformation of said constants.

The ESMA algorithm was designed with the hopes that the conjectures it discovers will reveal more efficient expansions of mathematical constants and more generally reveal unknown underlying patterns. Looking forward, the ESMA algorithm represents a wider effort to automate and accelerate mathematical discovery through the utilization of computational power. For instance, similar to how analyzing the algorithm's abilities required developing novel methods for manipulating CFs (e.g., the Folding transform), the proofs for its conjectures may require developing novel mathematical tools, further accelerating mathematical discovery.

The wider vision behind our approach is to find different ways of converting each target constant to an integer sequence (not necessarily through a CF), and to find different algorithms that can identify patterns in integer sequences. This approach may help find more efficient and direct discovery of conjectures on constants with many unknown properties, such as $G$ or $\zeta(5)$, potentially revealing their underlying structure or even a fundamental connection between them. We hope that this work will help inspire others to apply the ESMA algorithm and to develop other algorithms for the discovery of novel formulas of mathematical constants and identities between formulas.

# References

[1] S. R. Finch, Mathematical Constants, Cambridge University press, Encyclopedia of mathematics and its applications 94, 1-41, (2003).

[2] G. E. Andrews, An Introduction to Ramanujan's "lost" Notebook, The American Mathematical Monthly 86, 89-108, (1979).

[3] A quote from Berndt found in D. Peterson, Raiders of the Lost Notebook, University of Illinois at Urbana-Champaign, ATLAS, LAS News Magazine &Laquo, (2006).

[4] J. A. Harvey, Ramanujan's influence on string theory, blackholes and moonshine, In Philosophical Transactions of the Royal Society A: Mathematical, Physical and Engineering Sciences 378, (2019).

[5] R. Apéry, "Irrationalité de ζ (2) et ζ (3)", Astérisque 61, 11-13, (1979).

[6] G. Raayoni, S. Gottlieb, Y. Manor, et al. Generating conjectures on fundamental constants with the Ramanujan Machine, Nature 590, 67–73, (2021).

[7] J. M. Laughlin, N. J. Wyshinski, Real numbers with polynomial continued fraction expansions, Acta Arithmetica 116, 63-79, (2005).

[8] E. R. Berlekamp, Non-binary BCH decoding, North Carolina State University. Dept. of Statistics, Institute of Statistics mimeo series 502, (1966).

[9] J. L. Massey, Shift-Register Synthesis and BCH Decoding, IEEE Transaction on Information Theory 15, 122-127, (1969).

[10] C. E. Shannon, The Mathematical Theory of Communication, Bell System Technical Journal 27, 379-423, (1948).

[11] C. D. Olds, The Simple Continued Fraction Expression of e, American Mathematical M Monthly, 77, 968-974

[12] J. A. Reeds, N. J. A. Sloane, Shift-Register Synthesis (Modulo m)*, SIAM Journal on Computing 14, 505-513, (1985).

[13] A. Klein, Stream Ciphers, Springer London, (2013).

[14] G. H. Hardy, E. M. Wright, An Introduction to the Theory of Numbers 5 edition, Oxford University Press, (1980).

[15] G. Chrystal, Algebra: An Elementary Textbook vol. 2, Chelsea Publisher, (1964).

[16] H. S. Wall, Analytic Theory of Continued Fractions, D. Van Nostrand, New York: Chelsea, (1948).

[17] C. D. Olds, Continued Fractions, The Mathematical Association of America, (1963).

[18] F. G. Gustavson, Analysis of the Berlekamp-Massey Linear Feedback Shift-Register Synthesis Algorithm, IBM Journal of Research and Development, 20, 204-212, (1976).

[19] N. B. David, G. Nimri, On the Connection Between Irrationality Measures and Polynomial Continued Fractions, arXiv:2111.04468, (2021).

[20] D. Bowman, J. M. Laughlin, Polynomial Continued Fractions, Acta Arithmetica 103, 329-342, (2002).

[21] Y. Bugeaud, Distribution Modulo One and Diophantine Approximation (Cambridge Tracts in Mathematics), Cambridge University Press, (2012).

[22] S. Lang, Introduction to Diophantine Approximations, 2nd edition. Springer, New York, (1995).

[23] A. A. Cuyt, V. Petersen, B. Verdonk, H. Waadeland, W.B. Jones, Handbook of Continued Fractions for Special Functions, Springer Science & Business Media, (2008).

[24] L. J. Lange. An Elegant Continued Fraction for π, The American Mathematical Monthly, 106, 456, (1999).

[25] A. Denjoy, Sur une fonction réelle de Minkowski, Journal de Mathematiques pures et Appliquees 17, 105-155, (1938).

[26] L. Balkova, A. Hrušková, Continued Fractions of Square Roots of Natural Numbers, Acta Polytechnica, 53 (4), 322-328 (2013).

[27] https://mathworld.wolfram.com/BesselFunctionoftheFirstKind.html

# Appendix

## Appendix A: Additional Results by the ESMA Algorithm

In this section we present a small sample of conjectures and known formulas found by the ESMA algorithm. For each conjecture found, we present the signed $b_j$ sequence use for its extraction, the found LFSR or recurrence relation found for the $a_j$ sequence, the initial conditions of the recurrence relation, and the resultant $a_j$ sequence along with the CF's convergence rate. The convergence rate is a linear approximation of the number of digits obtained in the CF approximation per CF term in the log scale of approximation error. While super-exponential convergence is non-linear (in the log scale), this approximation captures the overall magnitude of convergence indicating what conjectures converge faster. Additional analysis of the convergence rate is best done by plotting the approximation error with CF depth.

As seen in the results, we solely obtained linear interlaced sequences in the partial denominator and all the results of constants found have polynomial CF expansions that converge super-exponentially. ESMA can additionally trivially catch second-degree algebraic numbers whose CF expansions converge exponentially.

We denote the initial conditions of the recurrence relation in list form under the relation. For example, if $a_0 = 1, a_1 = 1, a_2 = 1, a_3 = 3$ then we write the following under the recurrence relation: [1,1,1,3]. We further denote the periodic $b_j$ sequence using the following notation: $\{1, -1, \dots, 1\}$, this notation describes a period of the sequence thus fully describing the sequence, $b_{(n-1)\times\beta+1} = 1, b_{(n-1)\times\beta+2} = -1, \dots, \forall n \in \mathbb{N}\backslash 0$.

For our code refer to https://github.com/RamanujanMachine/RamanujanMachine/tree/master/ESMA. Notice, the code returns a sequence starting from $a_0$ while in our analysis we ignore $a_0$, moving it to the left-hand side of the equation.

| Novelty | Formula | $b_j$ | LFSR | $a_j$ | Convergence $\left[\frac{\text{digits}}{\text{term}}\right]$ |
|---|---|---|---|---|---|
| known | $-1+\mathrm{e} = 1+\cfrac{1}{1+\cfrac{1}{2+\cfrac{1}{1+\cdots}}}$ | $\{1,1,1\}$ | $a_j - 2a_{j-3} + a_{j-6} = 0$<br>[1,1,2,1,1,4] | $a_j = \begin{cases} 1 & j=3k \\ 1 & j=3k+1 \\ 2k & j=3k+2 \end{cases}$<br>$k \in \{0,1,2,\dots\}$ | 1.2868 |
| new and unproven | $\frac{1+\mathrm{e}}{4(-1+\mathrm{e})} = 1-\cfrac{1}{2+\cfrac{1}{5+\cfrac{1}{2-\cdots}}}$ | $\{-1,1,1,-1,1,1\}$ | $a_j - 2a_{j-6} + a_{j-12} = 0$<br>[1, 2, 5, 2, 3, 56, 5, 2, 21, 2, 7, 120] | $a_j = \begin{cases} 1+4k & j=6k \\ 2 & j=6k+1 \\ 5+16k & j=6k+2 \\ 2 & j=6k+3 \\ 3+4k & j=6k+4 \\ 56+64k & j=6k+5 \end{cases}$<br>$k \in \{0,1,2,\dots\}$ | 3.0142 |
| known | $\frac{1+\mathrm{e}}{-1+\mathrm{e}} = 2+\cfrac{1}{6+\cfrac{1}{10+\cfrac{1}{14+\cdots}}}$ | $\{1\}$ | $a_j - 2a_{j-1} + a_{j-2} = 0$<br>[2,6] | $a_j = 2+4j$ | 4.8512 |
| new and unproven | $\frac{-5+3\mathrm{e}}{3-\mathrm{e}} = 12-\cfrac{1}{1+\cfrac{1}{5-\cfrac{1}{1+\cdots}}}$ | $\{-1,1,-1,1\}$ | $a_j - 2a_{j-4} + a_{j-8} = 0$<br>[12, 1, 5, 1, 28, 1, 9, 1] | $a_j = \begin{cases} 12+16k & j=4k \\ 1 & j=4k+1 \\ 5+4k & j=4k+2 \\ 1 & j=4k+3 \end{cases}$<br>$k \in \{0,1,2,\dots\}$ | 2.1508 |
| new and unproven | $\frac{2+2\mathrm{e}}{-1+3\mathrm{e}} = 2-\cfrac{1}{1+\cfrac{1}{24+\cfrac{1}{3-\cdots}}}$ | $\{-1,1,1,-1,1,1\}$ | $a_j - 2a_{j-6} + a_{j-12} = 0$<br>[2, 1, 24, 3, 2, 13, 2, 5, 88, 7, 2, 29] | $a_j = \begin{cases} 2 & j=6k \\ 1+4k & j=6k+1 \\ 24+64k & j=6k+2 \\ 3+4k & j=6k+3 \\ 2 & j=6k+4 \\ 13+16k & j=6k+5 \end{cases}$<br>$k \in \{0,1,2,\dots\}$ | 3.4269 |
| new and unproven | $\frac{-3+5\mathrm{e}}{-6+6\mathrm{e}} = 1+\cfrac{1}{36+\cfrac{1}{2-\cfrac{1}{4-\cdots}}}$ | $\{1,1,-1,-1,1,1,-1,-1,1,1,-1,-1\}$ | $a_j - 2a_{j-12} + a_{j-24} = 0$<br>[2, 2, 18, 3, 3, 4, 3, 5, 2, 2, 6,2, 2, 6,54, 7, 3, 8, 3, 9, 2, 2, 10, 2] | $a_j = \begin{cases} 2 & j=12k \\ 2+4k & j=12k+1 \\ 18+36k & j=12k+2 \\ 3+4k & j=12k+3 \\ 3 & j=12k+4 \\ 4+4k & j=12k+5 \\ 3 & j=12k+6 \\ 5+4k & j=12k+7 \\ 2 & j=12k+8 \\ 2 & j=12k+9 \\ 6+4k & j=12k+10 \\ 2 & j=12k+11 \end{cases}$<br>$k \in \{0,1,2,\dots\}$ | 2.1106 |
| new and unproven | $\frac{1}{-2+2e^2} = 1 - \cfrac{1}{1 + \cfrac{1}{11 + \cfrac{1}{2 - \cdots}}}$ | $\{-1,1,1,-1,1,-1,1,1,-1,1\}$ | $a_j - 2a_{j-10} + a_{j-20} = 0$ <br> $[1, 1, 11, 2, 1, 3, 2, 36, 3, 4, 4,$ $5, 1, 3, 5, 84, 6, 4]$ | $a_j = \begin{cases} 1+3k & j=10k \\ 1 & j=10k+1 \\ 11+48k & j=10k+2 \\ 2+3k & j=10k+3 \\ 1 & j=10k+4 \\ 3 & j=10k+5 \\ 2+3k & j=10k+6 \\ 36+48k & j=10k+7 \\ 3+3k & j=10k+8 \\ 4 & j=10k+9 \end{cases}$ <br> $k \in \{0,1,2,\ldots\}$ | 2.3169 |
| known | $-\frac{1}{2} + \frac{e^2}{2} = 3 + \cfrac{1}{5 + \cfrac{1}{7 + \cfrac{1}{9 + \cdots}}}$ | $\{1\}$ | $a_j - 2a_{j-1} + a_{j-2} = 0$ <br> $[3,5]$ | $a_j = 3 + 2j$ | 4.2727 |
| new and unproven | $\frac{3+3\sqrt{e}}{2} = 3 + \cfrac{1}{2 - \cfrac{1}{1 + \cfrac{1}{36 - \cdots}}}$ | $\{1,-1,1,-1,1,-1\}$ | $a_j - 2a_{j-6} + a_{j-12} = 0$ <br> $[3,2,1,36,1,6,3,10,1,108,1,14]$ | $a_j = \begin{cases} 3 & j=6k \\ 2+8k & j=6k+1 \\ 1 & j=6k+2 \\ 36+72k & j=6k+3 \\ 1 & j=6k+4 \\ 6+8k & j=6k+5 \end{cases}$ <br> $k \in \{0,1,2,\ldots\}$ | 2.4322 |
| new and unproven | $\frac{2}{\tanh\left(\frac{1}{4}\right)} = 8 + \cfrac{1}{6 + \cfrac{1}{40 + \cfrac{1}{14 + \cdots}}}$ | $\{1,1\}$ | $a_j - 2a_{j-2} + a_{j-4} = 0$ <br> $[8,6,40,14]$ | $a_j = \begin{cases} 8+32k & j=2k \\ 6+8k & j=2k+1 \end{cases}$ <br> $k \in \{0,1,2,\ldots\}$ | 4.9003 |

**Table 1 | Conjectures involving $e, e^2, \sqrt{e}, tan(\frac{1}{4})$.** Digits per term value refers to the number of digits in the CF approximation per CF term, averaged over 100 terms.

| Novelty | Formula | $b_j$ | LFSR | $a_j$ | Convergence $\left[\frac{\text{digits}}{\text{term}}\right]$ |
|---|---|---|---|---|---|
| known | $\tan(1) = 1 + \cfrac{1}{1 + \cfrac{1}{1 + \cfrac{1}{3 + \cdots}}}$ | $\{1,1\}$ | $a_j - 2a_{j-2} + a_{j-4} = 0$ <br> $[1,1,1,3]$ | $a_j = \begin{cases} 1 & j=2k \\ 1+2k & j=2k+1 \end{cases}$ <br> $k \in \{0,1,2,\ldots\}$ | 1.8441 |
| known | $\frac{\tan(1)}{-1+\tan(1)} = 3 - \cfrac{1}{5 - \cfrac{1}{7 - \cfrac{1}{9 - \cdots}}}$ | $\{-1\}$ | $a_j - 2a_{j-1} + a_{j-2} = 0$ <br> $[3,5]$ | $a_j = 3 + 2j$ | 4.2727 |
| new and unproven | $\frac{2-\tan(1)}{-1+\tan(1)} = 1 - \cfrac{1}{4 + \cfrac{1}{2 - \cfrac{1}{1 + \cdots}}}$ | $\{-1,1,-1,1\}$ | $a_j - a_{j-1} + a_{j-2} - a_{j-3} - a_{j-4} + a_{j-5} - a_{j-6} + a_{j-7} = 0$ <br> $[1, 4, 2, 1, 5, 8, 2]$ | $a_j = \begin{cases} 1+k & j=4k \\ 3+k & j=4k+1 \\ 2 & j=4k+2 \\ 1 & j=4k+3 \end{cases}$ <br> $k \in \{0,1,2,\ldots\}$ | 1.8436 |
| new and unproven | $\frac{2}{\tan(1)} = 2 - \cfrac{1}{2 - \cfrac{1}{2 - \cfrac{1}{3 - \cdots}}}$ | $\{-1,-1,-1,-1\}$ | $a_j - 2a_{j-4} + a_{j-8} = 0$ <br> $[2, 2, 2, 3, 14, 5, 2, 6]$ | $a_j = \begin{cases} 2+12k & j=4k \\ 2+3k & j=4k+1 \\ 2 & j=4k+2 \\ 3+3k & j=4k+3 \end{cases}$ <br> $k \in \{0,1,2,\ldots\}$ | 3.0325 |
| new and unproven | $\frac{1}{-2+2\tan(1)} = 1 - \cfrac{1}{9 + \cfrac{1}{1 + \cfrac{1}{3 - \cdots}}}$ | $\{-1,1,1,-1,1,1\}$ | $a_j - 2a_{j-6} + a_{j-12} = 0$ <br> $[1, 9, 1, 3, 1, 1, 4, 21, 1, 6, 1, 1]$ | $a_j = \begin{cases} 1+3k & j=6k \\ 9+12k & j=6k+1 \\ 1 & j=6k+2 \\ 3+3k & j=6k+3 \\ 1 & j=6k+4 \\ 1 & j=6k+5 \end{cases}$ <br> $k \in \{0,1,2,\ldots\}$ | 1.8622 |
| new and unproven | $\frac{-2+2\tan(1)}{-3+2\tan(1)} = 10 - \cfrac{1}{4 - \cfrac{1}{2 - \cfrac{1}{5 - \cdots}}}$ | $\{-1,-1,-1,-1\}$ | $a_j - 2a_{j-4} + a_{j-8} = 0$ <br> $[10, 4, 2, 5, 22, 7, 2, 8]$ | $a_j = \begin{cases} 10+12k & j=4k \\ 4+3k & j=4k+1 \\ 2 & j=4k+2 \\ 5+3k & j=4k+3 \end{cases}$ <br> $k \in \{0,1,2,\ldots\}$ | 3.0555 |
| new and unproven | $\frac{-5+4\tan(1)}{-7+5\tan(1)} = 2 - \cfrac{1}{3 - \cfrac{1}{2 - \cfrac{1}{2 - \cdots}}}$ | $\{-1,-1,-1,-1,-1,-1,-1,-1,-1$ $-1,-1,-1\}$ | $a_j - 2a_{j-12} + a_{j-24} = 0$ <br> $[2, 3, 2, 2, 3, 27, 4, 3, 5, 3, 6, 2, 2, 7, 2, 2, 7, 63, 8, 3, 9, 3, 10, 2]$ | $a_j = \begin{cases} 2 & j=12k \\ 3+4k & j=12k+1 \\ 2 & j=12k+2 \\ 2 & j=12k+3 \\ 3+4k & j=12k+4 \\ 27+36k & j=12k+5 \\ 4+4k & j=12k+6 \\ 3 & j=12k+7 \\ 5+4k & j=12k+8 \\ 3 & j=12k+9 \\ 6+4k & j=12k+10 \\ 2 & j=12k+11 \end{cases}$ <br> $k \in \{0,1,2,\ldots\}$ | 1.844 |
| new and unproven | $\frac{-1+\tan(1)}{5-3\tan(1)} = 1 + \cfrac{1}{1 + \cfrac{1}{2 + \cfrac{1}{2 + \cdots}}}$ | $\{1,1,1,1,1,1\}$ | $a_j - 2a_{j-6} + a_{j-12} = 0$ <br> $[1, 1, 2, 2, 1, 16, 1, 4, 2, 5, 1, 28]$ | $a_j = \begin{cases} 1 & j=6k \\ 1+3k & j=6k+1 \\ 2 & j=6k+2 \\ 2+3k & j=6k+3 \\ 1 & j=6k+4 \\ 16+12k & j=6k+5 \end{cases}$ <br> $k \in \{0,1,2,\ldots\}$ | 1.8491 |
| new and unproven | $\frac{4-2\tan(1)}{7-4\tan(1)} = 2-\cfrac{1}{1+\cfrac{1}{6-\cfrac{1}{4-\cdots}}}$ | $\{-1,1,-1,-1,-1,1,-1,-1\}$ | $a_j - 2a_{j-8} + a_{j-16} = 0$ [2, 1, 6, 4, 2, 4, 2, 2, 2, 1, 18, 7, 2, 7, 2, 2] | $a_j = \begin{cases} 2 & j=8k \\ 1 & j=8k+1 \\ 6+12k & j=8k+2 \\ 4+3k & j=8k+3 \\ 2 & j=8k+4 \\ 4+3k & j=8k+5 \\ 2 & j=8k+6 \\ 2 & j=8k+7 \end{cases}$ $k \in \{0,1,2,\ldots\}$ | 1.3186 |
| new and unproven | $\frac{8-5\tan(1)}{-3+2\tan(1)} = 2-\cfrac{1}{7-\cfrac{1}{8+\cfrac{1}{2-\cdots}}}$ | $\{-1,-1,1,-1,-1,1\}$ | $a_j - a_{j-1} + a_{j-3} - a_{j-4} - a_{j-6} + a_{j-7} - a_{j-9} + a_{j+10} = 0$ [2, 7, 8, 2, 2, 1, 8, 13, 14, 2] | $a_j = \begin{cases} 2+6k & j=6k \\ 7+6k & j=6k+1 \\ 8+6k & j=6k+2 \\ 2 & j=6k+3 \\ 1 & j=6k+4 \\ 1 & j=6k+5 \end{cases}$ $k \in \{0,1,2,\ldots\}$ | 1.8851 |

**Table 2 | Conjectures for $tan(1)$ (in radians).**

| Novelty | Formula | $b_j$ | LFSR | $a_j$ | Convergence $\left[\frac{\text{digits}}{\text{term}}\right]$ |
|---|---|---|---|---|---|
| known | $\frac{J_0(1)}{J_1(1)} = 2-\cfrac{1}{4-\cfrac{1}{6-\cfrac{1}{8-\cdots}}}$ | $\{-1\}$ | $a_j - 2a_{j-1} - a_{j-2} = 0$ [2,4] | $a_j = 2+2j$ | 4.2669 |
| new and unproven | $-1+\frac{J_0(1)}{J_1(1)} = 1-\cfrac{1}{3+\cfrac{1}{1+\cfrac{1}{5-\cdots}}}$ | $\{-1,1,1\}$ | $a_j - 2a_{j-3} + a_{j-6} = 0$ [1, 3, 1, 5, 7, 1] | $a_j = \begin{cases} 1+4k & j=3k \\ 3+4k & j=3k+1 \\ 1 & j=3k+2 \end{cases}$ $k \in \{0,1,2,\ldots\}$ | 2.6247 |
| new and unproven | $\frac{\frac{-3J_1(1)}{2}+J_0(1)}{-J_0(1)+2J_1(1)} = 1-\cfrac{1}{11+\cfrac{1}{1+\cfrac{1}{3-\cdots}}}$ | $\{-1,1,1\}$ | $a_j - 2a_{j-3} - a_{j-6} = 0$ [1, 11, 1, 3, 19, 1] | $a_j = \begin{cases} 1+2k & j=3k \\ 11+7k & j=3k+1 \\ 1 & j=3k+2 \end{cases}$ $k \in \{0,1,2,\ldots\}$ | 2.6299 |
| new and unproven | $\frac{J_0(1)}{2J_1(1)} = 1-\cfrac{1}{8-\cfrac{1}{3-\cfrac{1}{16-\cdots}}}$ | $\{-1,-1\}$ | $a_j - 2a_{j-2} - a_{j-4} = 0$ [1,8,3,16] | $a_j = \begin{cases} 1+2k & j=2k \\ 8+16k & j=2k+1 \end{cases}$ $k \in \{0,1,2,\ldots\}$ | 4.2669 |
| new and unproven | $\frac{-J_0(1)+3J_1(1)}{-J_1(1)+J_0(1)} = 1+\cfrac{1}{2-\cfrac{1}{1+\cfrac{1}{1+\cdots}}}$ | $\{1,-1,1\}$ | $a_j - 2a_{j-3} - a_{j-6} = 0$ [1,2,1,1,3,1] | $a_j = \begin{cases} 1 & j=3k \\ 2+k & j=3k+1 \\ 1 & j=3k+2 \end{cases}$ $k \in \{0,1,2,\ldots\}$ | 1.1102 |
| new and unproven | $\frac{-2J_0(1)+4J_1(1)}{-4J_0(1)+7J_1(1)} = 12-\cfrac{1}{4-\cfrac{1}{20-\cfrac{1}{6-\cdots}}}$ | $\{-1\}$ | $a_j - 2a_{j-2} - a_{j-4} = 0$ [12,4,20,6] | $a_j = \begin{cases} 12+8k & j=2k \\ 4+2k & j=2k+1 \end{cases}$ $k \in \{0,1,2,\ldots\}$ | 4.2897 |
| new and unproven | $\frac{-12047J_1(1)+6928J_0(1)}{8(1777J_0(1)+3090J_1(1))} = 1-\cfrac{1}{112-\cfrac{1}{2-\cfrac{1}{143+\cdots}}}$ | $\{-1,-1,-1,1,1,-1,-1,-1,1,1,-1,-1,-1,1,1\}$ | $a_j - 2a_{j-15} - a_{j-30} = 0$ [1, 112, 2, 143, 1, 2, 2, 45, 1, 1, 3, 2, 53, 1, 1, 3, 240, 4, 271, 1, 4, 2, 77, 1, 1, 5, 2, 85, 1, 1] | $a_j = \begin{cases} 1+2k & j=15k \\ 112+128k & j=15k+1 \\ 2+2k & j=15k+2 \\ 143+128k & j=15k+3 \\ 1 & j=15k+4 \\ 2+2k & j=15k+5 \\ 2 & j=15k+6 \\ 45+32k & j=15k+7 \\ 1 & j=15k+8 \\ 1 & j=15k+9 \\ 3+2k & j=15k+10 \\ 2 & j=15k+11 \\ 53+32k & j=15k+12 \\ 1 & j=15k+13 \\ 1 & j=15k+14 \end{cases}$ $k \in \{0,1,2,\ldots\}$ | 2.021 |
| known | $\frac{J_1(1)}{J_2(1)} = 4-\cfrac{1}{6-\cfrac{1}{8-\cfrac{1}{10-\cdots}}}$ | $\{-1\}$ | $a_j - 2a_{j-1} - a_{j-2} = 0$ [4,6] | $a_j = 4+2j$ | 4.2784 |
| new and unproven | $\frac{J_0(1)}{J_2(1)} = 6+\cfrac{1}{1+\cfrac{1}{1+\cfrac{1}{1+\cdots}}}$ | $\{1,1,1,1\}$ | $a_j - 2a_{j-4} - a_{j-8} = 0$ [6,1,1,1,14,1,3,1] | $a_j = \begin{cases} 6+8k & j=4k \\ 1 & j=4k+1 \\ 1+2k & j=4k+2 \\ 1 & j=4k+3 \end{cases}$ $k \in \{0,1,2,\ldots\}$ | 1.8617 |
| new and unproven | $\frac{J_1(1)}{J_3(1)} = 23-\cfrac{1}{1+\cfrac{1}{1+\cfrac{1}{39-\cdots}}}$ | $\{-1,1,1\}$ | $a_j - 2a_{j-3} - a_{j-6} = 0$ [23,1,1,39,2,1] | $a_j = \begin{cases} 23+16k & j=3k \\ 1+k & j=3k+1 \\ 1 & j=3k+2 \end{cases}$ $k \in \{0,1,2,\ldots\}$ | 2.6693 |
| new and unproven | $\frac{J_3(1)}{J_5(1)} = 79-\cfrac{1}{1+\cfrac{1}{2+\cfrac{1}{27+\cdots}}}$ | $\{-1,1,1,1,1,-1,1,1,1,1,-1,1,1,1,1\}$ | $a_j - 2a_{j-15} + a_{j-30} = 0$ [79, 1, 2, 27, 2, 2, 1, 1, 35, 2, 2, 175, 1, 1, 1, 207, 3, 2, 59, 2, 4, 1, 1, 67, 2, 4, 303, 1, 3, 1] | $a_j = \begin{cases} 79+128k & j=15k \\ 1+2k & j=15k+1 \\ 2 & j=15k+2 \\ 27+32k & j=15k+3 \\ 2 & j=15k+4 \\ 2+2k & j=15k+5 \\ 1 & j=15k+6 \\ 1 & j=15k+7 \\ 35+32k & j=15k+8 \\ 2 & j=15k+9 \\ 2+2k & j=15k+10 \\ 175+128k & j=15k+11 \\ 1 & j=15k+12 \\ 1+2k & j=15k+13 \\ 1 & j=15k+14 \end{cases}$ $k \in \{0,1,2,\ldots\}$ | 2.0013 |
| new and unproven | $\frac{2}{\frac{J_5(1)}{J_3(1)}+1} = 2 - \frac{1}{40 - \frac{1}{3 - \frac{1}{56 - \cdots}}}$ | $\{-1,-1\}$ | $a_j - 2a_{j-2} + a_{j-4} = 0$ <br> [2,40,3,56] | $a_j = \begin{cases} 2+k & j=2k \\ 40+16k & j=2k+1 \end{cases}$ $k \in \{0,1,2,\ldots\}$ | 4.3008 |
| new and unproven | $\frac{1}{\frac{J_7(1)}{J_5(1)}+1} = 1 - \frac{1}{168 - \frac{1}{2 - \frac{1}{2 - \cdots}}}$ | $\{-1,-1,-1,$ $-1,-1,-1,$ $-1,-1,-1\}$ | $a_j - 2a_{j-9} + a_{j-18} = 0$ <br> [1, 168, 2, 2, 2, 25, 3, 2, 264, 2, 312,3, 2, 2, 41, 3, 3, 408] | $a_j = \begin{cases} 1+k & j=9k \\ 168+144k & j=9k+1 \\ 2+k & j=9k+2 \\ 2 & j=9k+3 \\ 2 & j=9k+4 \\ 25+16k & j=9k+5 \\ 3 & j=9k+6 \\ 2+k & j=9k+7 \\ 264+144k & j=9k+8 \end{cases}$ $k \in \{0,1,2,\ldots\}$ | 2.728 |

**Table 3 | Conjectures of ratios of Bessel functions of different order at $x = 1$** [27].

Overall, the ESMA algorithm was run on more than 50 constants, in the following list we present a sample of constants which the algorithm did not find conjectures for: $\zeta(2)$, $\zeta(3)$, $\zeta(5)$, $\pi$, $G$, $\sqrt{\varphi}$, $2^{\frac{1}{3}}$, $100^{\frac{1}{5}}$, and many more. We conjecture that these constants could not be caught by ESMA as they are not in the representable set of the SICF structure utilized (see section 4).

## Efficient Continued Fraction Expansions found by the ESMA Algorithm

The ESMA algorithm found various conjectures for constants that converge significantly faster than known formulas. Even for $2^{\text{nd}}$ degree algebraic numbers, whose CF can be trivially calculated, the ESMA algorithm was able to find more efficient expansions of various constants such as $\varphi, \sqrt{2}$, and many more. Below we denote a sample of conjectures found on $J_1(1)/J_3(1)$, $J_5(1)/J_3(1)$, and $\varphi$ which converge at a faster rate, and thus more efficiently approximate the constant, than their simple CF expansions (to our knowledge) [17].

Looking at $J_1(1)/J_3(1)$, it's simple CF expansion converges at a rate of 1.8643 digits per term whilst our conjecture on $J_1(1)/J_3(1)$ (see Table 3) converges at a rate of 2.6693 digits per term (averaged over 100 terms). In a similar manner, our conjecture on $2/((J_5(1)/J_3(1)) + 1)$ (see Table 3) converges at a rate of 4.3008 digits per term, whilst the simple CF for $J_5(1)/J_3(1)$ converges at a rate of 1.9046, more than 2 times as fast. Furthermore, the simple CF of The Golden ratio, $\varphi$, which is famously given by $\varphi = 1 + 1/(1 + 1/(1 + 1/(\ldots)))$ converges at a rate of 0.4096 digits per term, while the following conjecture found by ESMA converges at a rate of 2.4576 digits per term:

$$\frac{1+2\varphi}{-3+2\varphi} = 18 - \frac{1}{18 - \frac{1}{18 - \cdots}}. \quad (18)$$

The convergence plots for all the above conjectures and formulas are presented below.

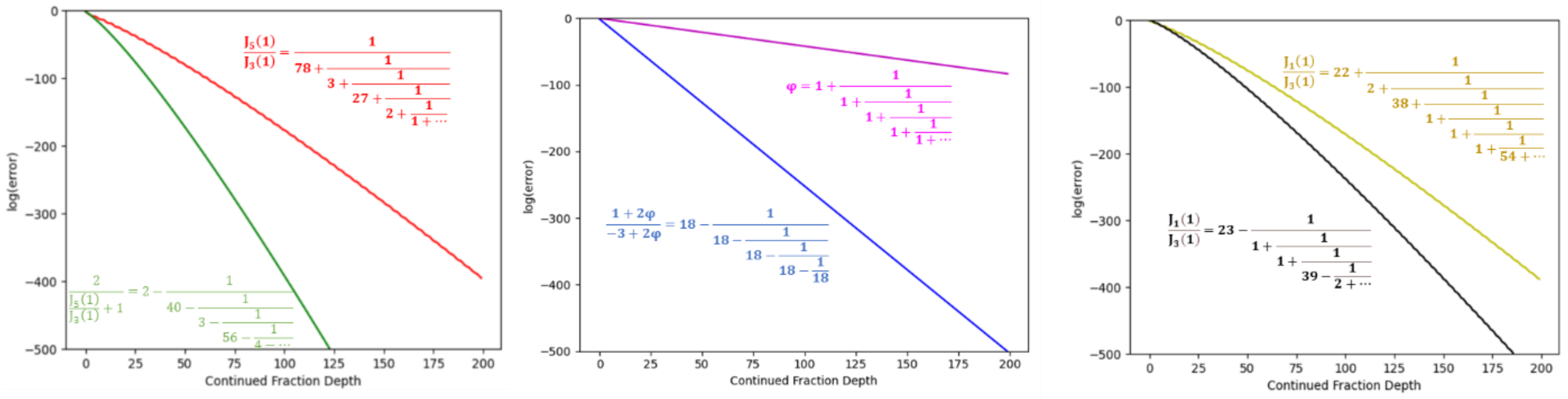

**Figure 4 | ESMA results that converges faster than the corresponding simple CFs**. log (error) plot of a known formula and conjecture on $\varphi$, $J_1(1)/J_3(1)$, and $J_5(1)/J_3(1)$. We see the conjectures found by the ESMA algorithm (green, blue, and black) approximate the constants significantly faster.

The increased rate of convergence of our conjectures relative to their respective simple CF formulas can be clearly seen in figure 4. We obtain conjectures for which less CF terms are needed for the same digit accuracy in approximating the constant.

All the conjectures are simply a Mobius transform of the constant we attempt to approximate. In these cases, the Mobius transforms do not "bother" us as we can find an expression which converges at the same rate while isolating the constant. For example, the conjecture on $\varphi$ can be written in the following way to directly approximate $\varphi$:

$$\varphi = -\frac{3}{2} + \cfrac{2}{17 - \cfrac{1}{18 - \cfrac{1}{18 - \cfrac{1}{18 - \cdots}}}}. \tag{19}$$

We note that some fast-converging expansions found by ESMA are known but were previously unappreciated for their fast convergence. For example, the formula found for $-1/2 + \mathrm{e}^2/2$ (See Table 1), which can be easily derived from the known CF of $\tanh(1)$, converges at a rate of 4.2727 digits per term whilst its known simple CF formula converges at a rate of 2.2806 digits per term. The known formulas found further support the validity of our method and adds confidence in our conjectures.

## Appendix B.1: Definitions

We recall, any matrix $M = \begin{pmatrix} a & b \\ c & d \end{pmatrix}$ with a non-zero determinant, $\det(M) = ad - bc \neq 0$, can represent the following Mobius transform:

$$M(x) = \begin{pmatrix} a & b \\ c & d \end{pmatrix}(x) = \frac{ax+b}{cx+d} \quad a, b, c, d \in \mathbb{Z}. \tag{20}$$

As such, each single fraction of a CF, e.g., $\frac{b}{a+\cdots}$, can be represented as a Mobius transform of the following form:

$$\begin{pmatrix} 0 & b \\ 1 & a \end{pmatrix}(x) = \frac{b}{a+x}. \tag{21}$$

We refer to a single fraction of a CF, (23), as a **layer** of the CF.

As such, any CF can be represented as a composition or product of Mobius transforms which act on 0, for example:

$$\frac{b_1}{a_1 + \frac{b_2}{a_2}} = \begin{pmatrix} 0 & b_1 \\ 1 & a_1 \end{pmatrix}\begin{pmatrix} 0 & b_2 \\ 1 & a_2 \end{pmatrix}(0).$$

Therefore, we can represent a CF as a matrix product of its layers in the following form:

$$\alpha = a_0 + \lim_{N\to\infty} \prod_{j=1}^{N} \begin{pmatrix} 0 & b_j \\ 1 & a_j \end{pmatrix}(0), \quad a_j, b_j \in \mathbb{Z}. \tag{22}$$

In this paper, an **Interlaced Continued Fraction** refers to a CF in which its partial numerator and denominator are interlaced sequences: a sequence formed of $\beta \geq 1$ different integer polynomial sub-sequences that alternate with a given order. An interlaced CF of a number $\alpha \in \mathbb{R}$ is a CF of the form:

$$\alpha = a_0 + \cfrac{b_1}{a_1 + \cfrac{b_2}{a_2 + \cfrac{b_3}{\ddots + \cfrac{b_{\beta+1}}{a_{\beta+1} + \cdots}}}} = a_0 + \cfrac{B_1(1)}{A_1(1) + \cfrac{B_2(1)}{A_2(1) + \cfrac{B_3(1)}{\ddots + \cfrac{B_1(2)}{A_1(2) + \cdots}}}} \tag{23}$$

Where $b_{(n-1)\times\beta+i} = B_i(n), a_{(n-1)\times\beta+i} = A_i(n),\ A_i, B_i \in \mathbb{Z}[x], \forall n \in \mathbb{N}\backslash 0$ are the partial numerator and denominator sequences, respectively, and $i \in \{1, .., \beta\}$ such that $\beta \in \mathbb{N}$ represents the amount of sub-sequences that are interlaced in the partial numerator/denominator, also referred to as the **period**. For each $i$, $A_i(n)$ represents the $n^{th}$ element of the $i^{th}$ polynomial sub-sequence. For the rest of this paper the leading integer $a_0$ will be ignored as it is more comfortable to treat it as part of the constant $\alpha$.

The $n^{th}$ **collapsed matrix** of an interlaced CF refers to a matrix, denoted by $M_n$, which is the product of all matrices in the $n^{th}$ period of the interlaced CF. Meaning given the above interlaced CF of $\alpha$ and of period $\beta$:

$$\alpha = \lim_{N\to\infty} \prod_{n=1}^{N} \prod_{i=1}^{\beta} \begin{pmatrix} 0 & b_{(n-1)\times\beta+i} \\ 1 & a_{(n-1)\times\beta+i} \end{pmatrix}(0) = \lim_{N\to\infty} \prod_{n=1}^{N} \prod_{i=1}^{\beta} \begin{pmatrix} 0 & B_i(n) \\ 1 & A_i(n) \end{pmatrix}(0) = \lim_{N\to\infty} \prod_{n=1}^{N} M_n(0) \tag{24}$$

$$M_n = \prod_{i=1}^{\beta} \begin{pmatrix} 0 & B_i(n) \\ 1 & A_i(n) \end{pmatrix} = \begin{pmatrix} c_n & d_n \\ e_n & f_n \end{pmatrix} \quad c, d, e, f \in \mathbb{Z}[x] \tag{25}$$

where $b_j$ and $a_j$ are the interlaced sequences from (23). All layers in each period are polynomial with the same $n$, as can be seen from the notation in (23), $b_{(n-1)\times\beta+i} = B_i(n), a_{(n-1)\times\beta+i} = A_i(n) \in \mathbb{Z}[x],\ n \in \mathbb{N}\backslash 0$. Therefore, the elements of the collapsed matrix $c, d, e, f$, being a product of various integer polynomial sequences with index $n$, will be polynomial with index $n$, and the collapsed matrix is a polynomial matrix.

## Appendix B.2: The Folding Transform and Proof of Theorem 1

We introduce the Folding Transform; a transform which acts on a constant's interlaced CF expansion, or in a more general sense acts on polynomial matrices. Given a constant $\alpha \in \mathbb{R}$ which has an interlaced CF expansion with a known period, we can "fold" (multiply) its layers in each period to form a polynomial matrix. Therefore, our analysis on general polynomial matrices applies to interlaced CFs in particular. The Folding transform applies a Mobius transform on $\alpha$ on the left-hand side of the equation while multiplying each of $\alpha$'s polynomial matrices on the right-hand side. The result of the Folding transform is a polynomial CF equal to a Mobius transform of $\alpha$.

Given an $\alpha$ which satisfies $\alpha = \lim_{N\to\infty}[\prod_{n=1}^{N} M_n](\infty)$, with a collapsed matrix $M_n = \begin{pmatrix} c_n & d_n \\ e_n & f_n \end{pmatrix}$ where $e_n \neq 0\ \forall n \in \mathbb{N}\backslash 0$, the Folding transform utilizes a sequence of matrices $U_n$ such that:

$$\alpha' = (M_1U_2)^{-1}(\alpha) = \lim_{N\to\infty}[\prod_{n=2}^{N} U_n^{-1}M_nU_{n+1}](\infty) \tag{26}$$

where $(M_1U_2)^{-1}(\alpha)$ is a Mobius transform acting on $\alpha$ and $U_n^{-1}M_nU_{n+1}$ is a polynomial CF layer. Two constants or CFs that are connected by the Folding transform are said to be **semi-equivalent**. This intriguing transform reveals a novel connection between polynomial CFs and interlaced CFs. Every constant that has an interlaced CF expansion, up to some Mobius transform, must have a polynomial CF expansion regardless of the number of interlaced sequences and the seeming complexity of the interlaced CF.

It is important to note that for every converging interlaced CFs as we defined them (see Appendix B.1), $e_n \neq 0$. Assessing the collapsed matrices, if we falsely assume $e_n = 0$ we obtain a CF which does not converge. If $M_n = \begin{pmatrix} c_n & d_n \\ 0 & f_n \end{pmatrix}\ \forall n \in \mathbb{N}\backslash 0$, the $n\beta - 1^{th}$ convergent (rational approximation of the CF) can be represented in the following form [6][23]:

$$\frac{p_{n\beta-1}}{q_{n\beta-1}} = \begin{pmatrix} p_{n\beta-1} & p_{n\beta} \\ q_{n\beta-1} & q_{n\beta} \end{pmatrix}(\infty) = \begin{pmatrix} 0 & b_1 \\ 1 & a_1 \end{pmatrix} \cdots \begin{pmatrix} 0 & b_{n\beta} \\ 1 & a_{n\beta} \end{pmatrix}(\infty) = M_1M_2 \dots M_n(\infty)$$

$$\begin{pmatrix} c_n & d_n \\ 0 & f_n \end{pmatrix}(\infty) = \infty$$

$$M_1M_2 \dots \begin{pmatrix} c_n & d_n \\ 0 & f_n \end{pmatrix}(\infty) = M_1M_2 \dots M_{n-1}(\infty) =^{\text{by induction}} \infty.$$

Proving the CF diverges. Therefore, any converging interlaced CF must satisfy $e_n \neq 0$ and one could apply the Folding transform on it. Note, we see here that for an interlaced CF, $[\prod_{n=1}^{N} M_n](\infty) = [\prod_{n=1}^{N-1} M_n](0)$, and therefore taking the mapping at $\infty$ is simply assessing a sub-sequence of the convergents the CF.

**Theorem 1.** Let $M_n = \begin{pmatrix} c_n & d_n \\ e_n & f_n \end{pmatrix}$, satisfying $e_n \neq 0$, be some polynomial matrix ($c, d, e, f \in \mathbb{Z}[x]$) for which the following limit exists $\lim_{N\to\infty}[\prod_{n=1}^{N} M_n](\infty)$. Then there exists some Mobius transform $T$ with entries in $\mathbb{Z}$, and polynomials $a', b' \in \mathbb{Z}[x]$ such that:

$$\lim_{N\to\infty}\left[T\prod_{n=1}^{N} M_n\right](\infty) = \mathbb{K}_1^{\infty}\frac{b'(n)}{a'(n)}$$

**Proof:**

For a given polynomial matrix $M_n$, satisfying $\lim_{N\to\infty}[\prod_{n=1}^{N} M_n](\infty) = \alpha$, we utilize the Folding transform to obtain a polynomial CF of a Mobius transform of $\alpha$ which we denote $\alpha'$. To prove Theorem 1, we must prove that:

1. The Folding transform does not change the limit of the polynomial matrix product: $\lim_{N\to\infty}[\prod_{n=1}^{N} M_n](\infty) = \alpha$.
2. The Folding transform obtains a polynomial CF.

1) To prove that the limit is unchanged we must show that just as the original product of polynomial matrices converges to $\alpha$ , $\lim_{N\to\infty}[\prod_{n=1}^{N} M_n](\infty) = \alpha$, so too does the resultant polynomial CF meaning we must prove that (see equation 26):

$$\lim_{N\to\infty}\left[U_1\prod_{n=1}^{N} U_n^{-1}M_nU_{n+1}\right](\infty) = \lim_{N\to\infty} U_1U_1^{-1}M_1U_2U_2^{-1}M_2 \dots M_NU_{N+1}(\infty) = \lim_{N\to\infty} M_1M_2 \dots M_NU_{N+1}(\infty) = \alpha$$

Given a polynomial matrix, $M_n = \begin{pmatrix} c_n & d_n \\ e_n & f_n \end{pmatrix}, c, d, e, f \in \mathbb{Z}[x], e_n \neq 0\ \forall n \in \mathbb{N}\backslash 0$, the Folding transform utilizes the following polynomial matrices $U_n$:

$$U_n = \begin{pmatrix} 1 & c_n \\ 0 & e_n \end{pmatrix},\ U_n^{-1} = \begin{pmatrix} e_n & -c_n \\ 0 & 1 \end{pmatrix}. \quad (27)$$

Assessing the result of mapping $U_{n+1}$ to $\infty$ we obtain, $\begin{pmatrix} 1 & c_{n+1} \\ 0 & e_{n+1} \end{pmatrix}(\infty) = \infty$.

Therefore, we find that the limit is unchanged by the Folding transform as

$$\lim_{N\to\infty} M_1 M_2 \dots M_N U_{N+1}(\infty) = \lim_{N\to\infty} M_1 M_2 \dots M_N(\infty) = \lim_{N\to\infty}\left[\prod_{n=1}^{N} M_n\right](\infty) = \alpha$$

2) We prove that we obtain a polynomial CF. Applying the Folding transform, the following product is obtained:

$$\alpha = \lim_{N\to\infty}\left[U_1 \prod_{n=1}^{N} U_n^{-1} M_n U_{n+1}\right](\infty) = \lim_{N\to\infty}\left[U_1 \prod_{n=1}^{N} \begin{pmatrix} e_n & -c_n \\ 0 & 1 \end{pmatrix}\begin{pmatrix} c_n & d_n \\ e_n & f_n \end{pmatrix}\begin{pmatrix} 1 & c_{n+1} \\ 0 & e_{n+1} \end{pmatrix}\right](\infty)$$

$$\alpha = \lim_{N\to\infty}\left[U_1 \prod_{n=1}^{N} \begin{pmatrix} 0 & e_{n+1}(-\det(M_n)) \\ e_n & e_n c_{n+1} + f_n e_{n+1} \end{pmatrix}\right](\infty).$$

We apply the equivalence transform [22-23] with the following sequence $g_n = \begin{cases} 1, n = 0 \\ e_n, n > 0 \end{cases}$ to normalize the bottom left element of each matrix in the product to obtain correct CF form (see equation (21) in Appendix B.1):

$$\alpha = \lim_{N\to\infty}\left[U_1 \prod_{n=1}^{N} \begin{pmatrix} 0 & g_{n-1} g_n e_{n+1}(e_n d_n - c_n f_n) \\ e_n & g_n(e_n c_{n+1} + f_n e_{n+1}) \end{pmatrix}\right](\infty) = \lim_{N\to\infty}\left[U_1 \prod_{n=1}^{N} \begin{pmatrix} 0 & e_{n-1} e_{n+1}(e_n d_n - c_n f_n) \\ 1 & (e_n c_{n+1} + f_n e_{n+1}) \end{pmatrix}\right](\infty)$$

Note, The equivalence transform [23] refers to the following transform: Given a CF $\{a_j, b_j\}$ and any non-zero infinite sequence of $\{g_i\}_{i=1}^{\infty} \in \mathbb{C}$:

$$a_0 + \frac{b_1}{a_1 + \frac{b_2}{a_2 + \cdots}} = a_0 + \frac{g_1 b_1}{g_1 a_1 + \frac{g_1 g_2 b_2}{g_2 a_2 + \cdots}},$$

These two CFs are equivalent to any depth $n \in \mathbb{N}$.

Since $c, d, e, f \in \mathbb{Z}[x]$, are all polynomial with $n$ any product or sum of the sequences is polynomial with $n$. The first layer of $n = 1$ is taken out of the product as it doesn't follow the polynomial pattern of the rest of the CF, as $g_0 = 1$, and we obtain:

$$\alpha = \lim_{N\to\infty}\left[\begin{pmatrix} 1 & c_1 \\ 0 & e_1 \end{pmatrix}\begin{pmatrix} 0 & e_2(e_1 d_1 - c_1 f_1) \\ 1 & (e_1 c_2 + f_1 e_2) \end{pmatrix} \prod_{n=2}^{N} \begin{pmatrix} 0 & e_{n-1} e_{n+1}(e_n d_n - c_n f_n) \\ 1 & (e_n c_{n+1} + f_n e_{n+1}) \end{pmatrix}\right](\infty).$$

Where, $\begin{pmatrix} 1 & c_1 \\ 0 & e_1 \end{pmatrix}\begin{pmatrix} 0 & e_2(e_1 d_1 - c_1 f_1) \\ 1 & (e_1 c_2 + f_1 e_2) \end{pmatrix} = U_1 U_1^{-1} M_1 U_2 = M_1 U_2$. We obtain a polynomial CF:

$$(M_1 U_2)^{-1}(\alpha) = \lim_{N\to\infty}\left[\prod_{n=2}^{N} \begin{pmatrix} 0 & e_{n-1} e_{n+1}(e_n d_n - c_n f_n) \\ 1 & (e_n c_{n+1} + f_n e_{n+1}) \end{pmatrix}\right](\infty). \quad (28)$$

The above polynomial CF is said to be semi-equivalent to $\alpha$'s interlaced CF.

We can utilize the connection revealed by the Folding transform to assess properties of interlaced CFs based on their semi-equivalent and more extensively researched polynomial CFs expansions. For example, given a polynomial CF that converges to an irrational limit $(M_1 U_2)^{-1}(\alpha)$ and that it is semi-equivalent to an interlaced CF that converges to $\alpha$, we can trivially deduce that $\alpha$ must be irrational.

## Appendix B.3: Applying the Folding Transform on the Simple CF of e

We apply the Folding transform on the simple CF of e (found in [11]). Looking at the simple CF of e, we can represent it in matrix form and notice that the CF is a simple interlaced CF (see equation (1)):

$$\mathrm{e}-2=\cfrac{1}{1+\cfrac{1}{2+\cfrac{1}{1+\cfrac{1}{1+\cfrac{1}{4+\cfrac{1}{1+\cdots}}}}}}.$$

$$\mathrm{e}-2=\left[\begin{pmatrix}0&1\\1&1\end{pmatrix}\begin{pmatrix}0&1\\1&2\end{pmatrix}\begin{pmatrix}0&1\\1&1\end{pmatrix}\begin{pmatrix}0&1\\1&1\end{pmatrix}\begin{pmatrix}0&1\\1&4\end{pmatrix}\begin{pmatrix}0&1\\1&1\end{pmatrix}\begin{pmatrix}0&1\\1&1\end{pmatrix}\begin{pmatrix}0&1\\1&6\end{pmatrix}\begin{pmatrix}0&1\\1&1\end{pmatrix}\ldots\right](0)$$

The above simple interlaced CF has a period of 3 (as clarified by the highlighted colors) with $B_i(n)=1\ \forall n,i\in\mathbb{N}$ and $A_i(n)=\begin{cases}1 & i=1\\ 2n & i=2\\ 1 & i=3\end{cases}\ \forall n\in\mathbb{N}$. We find its collapsed matrices and apply the Folding transform:

$$\mathrm{e}-2=\lim_{N\to\infty}\left[\prod_{n=1}^{N}\prod_{i=1}^{3}\begin{pmatrix}0&B_i(n)\\1&A_i(n)\end{pmatrix}\right](0)=\lim_{N\to\infty}\left[\prod_{n=1}^{N}\begin{pmatrix}2n&2n+1\\2n+1&2n+2\end{pmatrix}\right](0).$$

$M_n=\begin{pmatrix}2n&2n+1\\2n+1&2n+2\end{pmatrix}$ and therefore from equation (27), $U_n=\begin{pmatrix}1&2n\\0&2n+1\end{pmatrix}$, and the Folding transform is given by:

$$(M_1U_2)^{-1}(\mathrm{e}-2)=\lim_{N\to\infty}\left[\prod_{n=1}^{N}\begin{pmatrix}2n+1&-(2n)\\0&1\end{pmatrix}\begin{pmatrix}2n&2n+1\\2n+1&2n+2\end{pmatrix}\begin{pmatrix}1&2n+2\\0&2n+3\end{pmatrix}\right](0),$$

$$\left(\begin{pmatrix}0&1\\1&2\end{pmatrix}\begin{pmatrix}1&2\\0&3\end{pmatrix}\right)^{-1}(\mathrm{e}-2)=\lim_{N\to\infty}\left[\prod_{n=2}^{N}\begin{pmatrix}0&2n+3\\2n+1&8n^2+16n+8\end{pmatrix}\right](0).$$

We then apply equivalence theorem (see note in Appendix B.1) to reach the correct form of a CF layer and we get the following polynomial CF:

$$\frac{-8\mathrm{e}+19}{\mathrm{e}-2}=\lim_{N\to\infty}\left[\prod_{n=2}^{N}\begin{pmatrix}0&4n^2+4n-3\\1&8n^2+16n+8\end{pmatrix}\right](0),$$

$$\frac{-8\mathrm{e}+19}{\mathrm{e}-2}=\cfrac{21}{72+\cfrac{45}{128+\cfrac{77}{200+\cdots}}}. \tag{29}$$

## Appendix B.4: Determinant Property of the Folding Transform

When applying the Folding transform one can notice the determinant property of the Folding transform. Notice that in the Folding transform's product $U_n^{-1}M_nU_{n+1}$, in the left-hand side product we obtain $U_n^{-1}M_n=\begin{pmatrix}e_n&-c_n\\0&1\end{pmatrix}\begin{pmatrix}c_n&d_n\\e_n&f_n\end{pmatrix}=\begin{pmatrix}0&-\det(M_n)\\e_n&f_n\end{pmatrix}$ where the top right element of the resultant matrix is equal to minus the determinant of the collapsed matrix: $-\left|\begin{pmatrix}c_n&d_n\\e_n&f_n\end{pmatrix}\right|=e_nd_n-c_nf_n$. An interlaced CF of period $\beta$'s collapsed matrix determinant at each period is given by:

$$\left|\begin{pmatrix}c_n&d_n\\e_n&f_n\end{pmatrix}\right|=\left|\begin{pmatrix}0&B_1(n)\\1&A_1(n)\end{pmatrix}\begin{pmatrix}0&B_2(n)\\1&A_2(n)\end{pmatrix}\cdots\begin{pmatrix}0&B_\beta(n)\\1&A_\beta(n)\end{pmatrix}\right|=\prod_{i=1}^{\beta}(-B_i(n))=(-1)^\beta\prod_{i=1}^{\beta}B_i(n).$$

Therefore, we can simplify the left-hand side product to:

$$\begin{pmatrix}e_n&-c_n\\0&1\end{pmatrix}\begin{pmatrix}c_n&d_n\\e_n&f_n\end{pmatrix}=\begin{pmatrix}0&(-1)^\beta\prod_{i=1}^{\beta}B_i(n)\\e_n&f_n\end{pmatrix}. \tag{30}$$

This property eases calculations when trying to predict resultant polynomial CF properties such as polynomial degree after applying the Folding transform.

## Appendix C.1: Proof of Theorem 2

Recall, a simple CF can be represented in the following forms [17]

$$\alpha=a_0+\cfrac{1}{a_1+\cfrac{1}{a_2+\cfrac{1}{\cdots}}}=[a_0;a_1,a_2,a_3,a_4\ldots]. \tag{31}$$

To represent the partial numerator, we add a sign variable to the above notation. This enables us to represent **Signed Interlaced Continued Fractions (SICFs):**

$$\alpha=a_0+\cfrac{b_1}{a_1+\cfrac{b_2}{a_2+\cfrac{b_3}{\cdots}}},\ a_j>0\in\mathbb{Z},\ b_j\in\{1,-1\}\,\forall j\in\mathbb{N}, \tag{32}$$

$$\alpha = [a_0; (a_1, b_1), (a_2, b_2), (a_3, b_3), \dots].$$

For example, given the following SICF:

$$\alpha = a_0 - \frac{1}{a_1 + \frac{1}{a_2 - \frac{1}{\dots}}},$$

We can represent $\alpha$ in reduced notation:

$$\alpha = [a_0; (a_1, -1), (a_2, 1), (a_3, -1), \dots].$$

For the sake of our analysis, as previously done, we will ignore $a_0$ and represent each SICF in the following form:

$$\alpha = [(a_1, b_1), (a_2, b_2), (a_3, b_3), \dots]. \quad (33)$$

**Theorem 2**. For any signed interlaced continued fraction equal to $\alpha \in \mathbb{R} \setminus \mathbb{Q}$, there exists some Mobius transform of $\alpha$ which is equal to a simple interlaced continued fraction.

**Proof:**

SICFs with constant partial denominator sequences

We show that given a converging SICF of period $\beta$ with constant partial denominator sequences, meaning $a_j = a_{j+\beta n}\ \forall j \in \mathbb{N}, \forall n \in \mathbb{N}$, we converge to $2^{\text{nd}}$ degree algebraic numbers and therefore by [17] there exists a Mobius transform which is equal to a simple interlaced CF (a simple CF with a pattern in its partial denominator sequence). A constant SICF of period $\beta$ has a constant collapsed matrix:

$$M_n = M = \begin{pmatrix} a & b \\ c & d \end{pmatrix}, \qquad a, b, c, d \in \mathbb{Z}\ \forall n \in \mathbb{N}$$

As any simple or signed CF layers are made of unimodular matrices, and the product of unimodular matrices is a unimodular matrix, the determinant of $M$ is given by $\det(M) = ad - bc = \pm 1$. We want to study the limit of $M^k(0)$ which is known to converge, $\alpha = \lim_{k\to\infty} p_k/q_k = \lim_{k\to\infty} M^k(0)$, therefore in the case that $\det(M) = -1$ we can instead consider the limit of the sub-sequence $M^{2k}(0)$ for which $\det(M^2) = 1$, allowing us to assume $\det(M) = 1$. We denote the eigen vectors of matrix M as $\mathbf{v_1} = \begin{pmatrix} v_{11} \\ v_{12} \end{pmatrix}$ and $\mathbf{v_2} = \begin{pmatrix} v_{21} \\ v_{22} \end{pmatrix}$and the eigen values as $\lambda_1$ and $\lambda_2$ which satisfy $\lambda_1 \lambda_2 = \det(M) = \pm 1$. The eigen values are the roots of the second-degree monic characteristic polynomials of $M$ and therefore are $2^{\text{nd}}$ degree algebraic numbers.

If $\lambda_1 = \overline{\lambda_2} \in \mathbb{C}$

If the eigen values are complex, one must be the complex conjugate of the other. Complex eigen values indicate rotation and scaling and therefore the matrix $M^k$ converges to a rational number or diverges. We show that any real 2 by 2 matrix $M$ with complex eigen values is similar to a shift rotation matrix $R = \begin{pmatrix} r\cos(\theta) & -r\sin(\theta) \\ r\sin(\theta) & r\cos(\theta) \end{pmatrix}$.

Utilizing a single eigen vector, $\mathbf{v_1}$, we construct a matrix $C = \begin{pmatrix} \mathrm{Re}(v_{11}) & -\mathrm{Im}(v_{11}) \\ \mathrm{Re}(v_{12}) & -\mathrm{Im}(v_{12}) \end{pmatrix} = (\mathrm{Re}(\mathbf{v_1}), -\mathrm{Im}(\mathbf{v_1}))$, where Re and Im are the real and imaginary parts of a complex number respectively. Through simple matrix multiplication one can see that $M = CRC^{-1}$ by showing that $MC = CR$.

We therefore obtain a rotation matrix $R$ which rotates vectors on the unit circle. Looking at the $k = n\beta$ convergent: $M^k(0) = CR^kC^{-1}(0)$.

For all $\theta = 0, 2\pi l\ \ \forall l \in \mathbb{N}$ we obtain an identity matrix and converge to rational numbers which are not of interest to us. Otherwise, $M^k$ infinitely rotates the $C^{-1}(0)$ vector by $\theta$ and the CF does not converge.

If $|\lambda_1| \neq |\lambda_2|,\ \lambda_1, \lambda_2 \in \mathbb{R}$

If matrix $M$ has two distinct eigen values, then it is diagonalizable and can be written in the following form:

$$\begin{pmatrix} a & b \\ c & d \end{pmatrix} = PDP^{-1} = \frac{1}{v_{11}v_{22} - v_{12}v_{21}} \begin{pmatrix} v_{11} & v_{21} \\ v_{12} & v_{22} \end{pmatrix} \begin{pmatrix} \lambda_1 & 0 \\ 0 & \lambda_2 \end{pmatrix} \begin{pmatrix} v_{22} & -v_{21} \\ -v_{12} & v_{11} \end{pmatrix}.$$

Note, the above matrix represents a Mobius transform and therefore, $cM = M$.

Without loss of generality, we can assume $|\lambda_1| > |\lambda_2|$. Looking at $M^k$:

$$M^k(0) = PD^kP^{-1}(0) = \begin{pmatrix} v_{11} & v_{21} \\ v_{12} & v_{22} \end{pmatrix}\begin{pmatrix} \lambda_1^k & 0 \\ 0 & \lambda_2^k \end{pmatrix}\left(-\frac{v_{21}}{v_{11}}\right) = \begin{pmatrix} \lambda_1^k v_{11} & \lambda_2^k v_{21} \\ \lambda_1^k v_{12} & \lambda_2^k v_{22} \end{pmatrix}\left(-\frac{v_{21}}{v_{11}}\right),$$

$$M^k(0) = \frac{\lambda_1^k v_{11}\left(-\frac{v_{21}}{v_{11}}\right)+\lambda_2^k v_{21}}{\lambda_1^k v_{12}\left(-\frac{v_{21}}{v_{11}}\right)+\lambda_2^k v_{22}} = \frac{v_{11}\left(-\frac{v_{21}}{v_{11}}\right)+\frac{\lambda_2^k}{\lambda_1^k}v_{21}}{v_{12}\left(-\frac{v_{21}}{v_{11}}\right)+\frac{\lambda_2^k}{\lambda_1^k}v_{22}} \xrightarrow{k\to\infty}= \frac{v_{11}}{v_{12}}.$$

The eigen vector must satisfy $\begin{pmatrix} a-\lambda_1 & b \\ c & d-\lambda_1 \end{pmatrix}\begin{pmatrix} v_{11} \\ v_{12} \end{pmatrix} = 0$, therefore, the limit is given by: $\frac{v_{11}}{v_{12}} = \frac{\lambda_1 - d}{c}$.

As such, we have shown that in this case the SICF converges to the eigen values which are $2^{\text{nd}}$ degree algebraic numbers.

<u>If $|\lambda_1| = |\lambda_2|$, $\lambda_1, \lambda_2 \in \mathbb{R}$</u>

In this case $M$ is either diagonalizable or non-diagonalizable. We denote $|\lambda_1| = |\lambda_2|$ , and recall that $\lambda_1\lambda_2 = \det(M) = 1$ meaning the eigen values are of the same sign, $\lambda_1 = \lambda_2 = \lambda = \pm 1$, without loss of generality we assume $\lambda = 1$ as if $\lambda = -1$ we can look at $M^2$. If $M$ is diagonalizable, then $M$ is a mapping that is equivalent to scaling and we are guaranteed to either diverge or converge to a rational number:

$$M^k = PD^kP^{-1} = P\begin{pmatrix} \lambda & 0 \\ 0 & \lambda \end{pmatrix}^k P^{-1} = P\begin{pmatrix} 1 & 0 \\ 0 & 1 \end{pmatrix}^k P^{-1} = PP^{-1} = I.$$

This is an identity matrix and is therefore of no interest to us in the analysis of CF of irrational numbers. If $M$ is not diagonalizable, it must be of the following form: $P\begin{pmatrix} 1 & b \\ 0 & 1 \end{pmatrix}P^{-1}$. This is the simple Jordan form. Where $P$'s columns are given by the eigen vectors of matrix $M$. If $M = P\begin{pmatrix} 1 & b \\ 0 & 1 \end{pmatrix}P^{-1}$ , we see that $\begin{pmatrix} 1 & b \\ 0 & 1 \end{pmatrix}^k(x) = x + kb$ and therefore we obtain a Mobius transform that maps any $x$ to either infinity or to a rational number when applied repeatedly,

$$\lim_{k\to\infty} PM^kP^{-1}(0) = \lim_{k\to\infty} P\left(-\frac{v_{21}}{v_{11}} + kb\right) = \lim_{k\to\infty}\begin{pmatrix} v_{11} & v_{21} \\ v_{12} & v_{22} \end{pmatrix}\left(-\frac{v_{21}}{v_{11}} + kb\right) = \lim_{k\to\infty}\frac{kb}{v_{12}kb + v_{22}}.$$

If $v_{12} \neq 0$ we will converge to a rational number otherwise we approach infinity.

Therefore, in the cases where the SICF with constant sequences converges to an irrational number, the SICF converges to a $2^{\text{nd}}$ degree algebraic number. This ensures there is a simple interlaced CF expansion for the number, where the interlaced sequences are constant sequences. Simple CFs with constant sequences (or k-periodic CFs) have a bijective correspondence to the real $2^{\text{nd}}$ degree algebraic numbers as given by Minkowski's Question Mark Function [25-26], up to initial CF layers which do not always follow the periodic pattern. Therefore, all $2^{\text{nd}}$ degree algebraic numbers, up to some Mobius transform, have a mathematical formula involving a k-periodic CF [17]. As a result, as the SICF converges to a $2^{\text{nd}}$ degree algebraic number we are guaranteed that this same number will have a simple interlaced CF expansion (up to some Mobius transform).

Note, these CFs with constant sequences converge at an exponential rate [6]. Conjectures on $2^{\text{nd}}$ degree algebraic numbers can be trivially calculated [26] and thus further analysis of periodic CFs does not provide additional insight on the algorithm.

<u>SICFs with non-constant partial denominator sequences</u>

Given an SICF which converges to $\alpha \in \mathbb{R}\backslash\mathbb{Q}$,

$$\alpha = [(a_1, b_1), (a_2, b_2), (a_3, b_3), \dots, (a_\beta, b_\beta), \dots],\ b_j \in \{\pm 1\},$$

We prove there exists some Mobius transform of $\alpha$ which equals to a simple interlaced CF. We assume $\exists i \in \{1, \beta\}$ s.t $\deg(A_i) > 0$ , $A_i(n) = a_{(n-1)\times\beta + i}$ , meaning we have non-constant sequences in the partial denominator. We address base cases and prove two lemmas to prove the above theorem. We utilize the following identities to aid our proof.

<u>**Identities for proof**</u>

For a given SICF, we may find a negative partial numerator in the CF, e.g., $b_{j+1} = -1$. Our goal is to find an equivalent CF expansion with all partial numerators all equal to 1. To this end we utilize the following matrix identity:

$$(*)\begin{pmatrix}0 & 1\\ 1 & a_j\end{pmatrix}\begin{pmatrix}0 & -1\\ 1 & a_{j+1}\end{pmatrix}=\begin{pmatrix}0 & 1\\ 1 & a_j-1\end{pmatrix}\begin{pmatrix}1 & 1\\ 0 & 1\end{pmatrix}\begin{pmatrix}-1 & 0\\ 1 & 1\end{pmatrix}\begin{pmatrix}0 & 1\\ 1 & a_{j+1}-1\end{pmatrix}=\begin{pmatrix}0 & 1\\ 1 & a_j-1\end{pmatrix}\begin{pmatrix}0 & 1\\ 1 & 1\end{pmatrix}\begin{pmatrix}0 & 1\\ 1 & a_{j+1}-1\end{pmatrix}$$

We thus found a method to "get rid" of negative partial numerators but we pay for it by decreasing the denominators. In a case where all interlaced sequences only have values greater than 2 then we can simply apply (*) iteratively over all the SICF to obtain a simple interlaced CF. Therefore, we must address cases were this is not met.

We introduce identities which will aid us in our proof. Most of the identities are simply a result of matrix multiplication and decomposition, with identity 1 generalizing the decomposition in (*). For example, the derivation of identity 3 is given by:

$$\begin{pmatrix}0 & b_{j-1}\\ 1 & a_{j-1}\end{pmatrix}\begin{pmatrix}0 & b_j\\ 1 & 0\end{pmatrix}\begin{pmatrix}0 & b_{j+1}\\ 1 & a_{j+1}\end{pmatrix}=\begin{pmatrix}0 & b_{j-1}b_{j+1}\\ b_j & b_j a_{j+1}+b_{j+1}a_{j-1}\end{pmatrix}\overset{b_j^2=1}{=}_{\text{matrix is a Mobius transform}}\begin{pmatrix}0 & b_{j-1}b_j b_{j+1}\\ 1 & a_{j+1}+b_j b_{j+1}a_{j-1}\end{pmatrix}.$$

## Identity 1

$$\alpha=[\dots,(a_{j-1},b_{j-1}),(a_j,-1),\dots]=[\dots,(a_{j-1}-1,b_{j-1}),(1,1),(a_j-1,1),\dots].$$

## Identity 2

$$\alpha=[\dots,(a_{j-1},b_{j-1}),(0,b_j),(a_{j+1},b_{j+1}),\dots]=[\dots,(a_{j+1}+b_{j+1}b_j a_{j-1},b_{j-1}b_j b_{j+1}),\dots].$$

## Identity 3 – Equivalence Theorem

For any non-zero infinite sequence, $g_i$, the following equality is met (see note Appendix B.1)

$$\alpha=[\dots,(a_{j-1},\ b_{j-1}),(a_j,\ b_j),(a_{j+1},\ b_{j+1}),\dots]\to$$

$$\alpha=[\dots,(g_{j-1}a_{j-1},g_{j-2}g_{j-1}b_{j-1}),(g_j a_j,\ g_{j-1}g_j b_j),(g_{j+1}a_{j+1},\ g_j g_{j+1}b_{j+1}),\dots]$$

When applying identity 1 and 2 we assume it is applied in every period of the SICF, we will prove this does not affect convergence in the proof of lemma 1. We further note that if we **"shift"** our period in an SICF, it is simply equivalent to applying a Mobius transform on our constant $\alpha$. To shift our period is to look at the period in a different order, for example:

$$\alpha=\lim_{N\to\infty}\left[\prod_{n=1}^{N}\begin{pmatrix}0 & b_1\\ 1 & A_1(n)\end{pmatrix}\cdots\begin{pmatrix}0 & b_\beta\\ 1 & A_\beta(n)\end{pmatrix}\right](0)\to$$

$$\left(\begin{pmatrix}0 & b_1\\ 1 & A_1(n+1)\end{pmatrix}\right)^{-1}(\alpha)=\lim_{N\to\infty}\left[\prod_{n=1}^{N}\begin{pmatrix}0 & b_2\\ 1 & A_2(n)\end{pmatrix}\cdots\begin{pmatrix}0 & b_\beta\\ 1 & A_\beta(n)\end{pmatrix}\begin{pmatrix}0 & b_1\\ 1 & A_1(n+1)\end{pmatrix}\right](0).$$

## Proof of Theorem 2 for the General Case

We prove inductively that any SICF can be represented as a simple interlaced CF up to a Mobius transform. This will follow immediately if we can prove the next lemma.

1. If all $a_j$'s are positive, $a_j>0$, **we can decrease the number of $b_j$'s$=-1$** in each period **at least by 1**.

To prove the above lemma, we need to prove an additional lemma which addresses edge cases met in the first lemma:

2. Given an SICF for which we allow a single $a_j$ in each period to be zero. If there is a single $a_k=0$ in each period ($a_k=a_{(n-1)\times\beta+k}=0\ \forall n\in\mathbb{N}\backslash 0$), satisfying $a_k=0\ b_k=\pm 1$, while all other $a_j$'s are positive: then there is an equivalent representation of the CF with all $a_j$'s being positive **without increasing the number of $b_j$'s$=-1$** in each period.

### Proof of Lemma 2

Given a SICF of period $\beta$,

$$\alpha=[(a_1,b_1),(a_2,b_2),\dots,(a_\beta,b_\beta),\dots],\ b_j\in\{\pm 1\},$$

for which there exists a $k \in \{1, \dots, \beta\}$ such that $a_{(n-1)\times\beta\,+\,k} = 0\ \ \forall n \in \mathbb{N}\backslash 0$ we find an equivalent representation of the CF with all $a_j$'s$> 0$ without increasing the number of $b_j$'s$= -1$ in each period.

The proof of lemma 2 is described in the following flow chart:

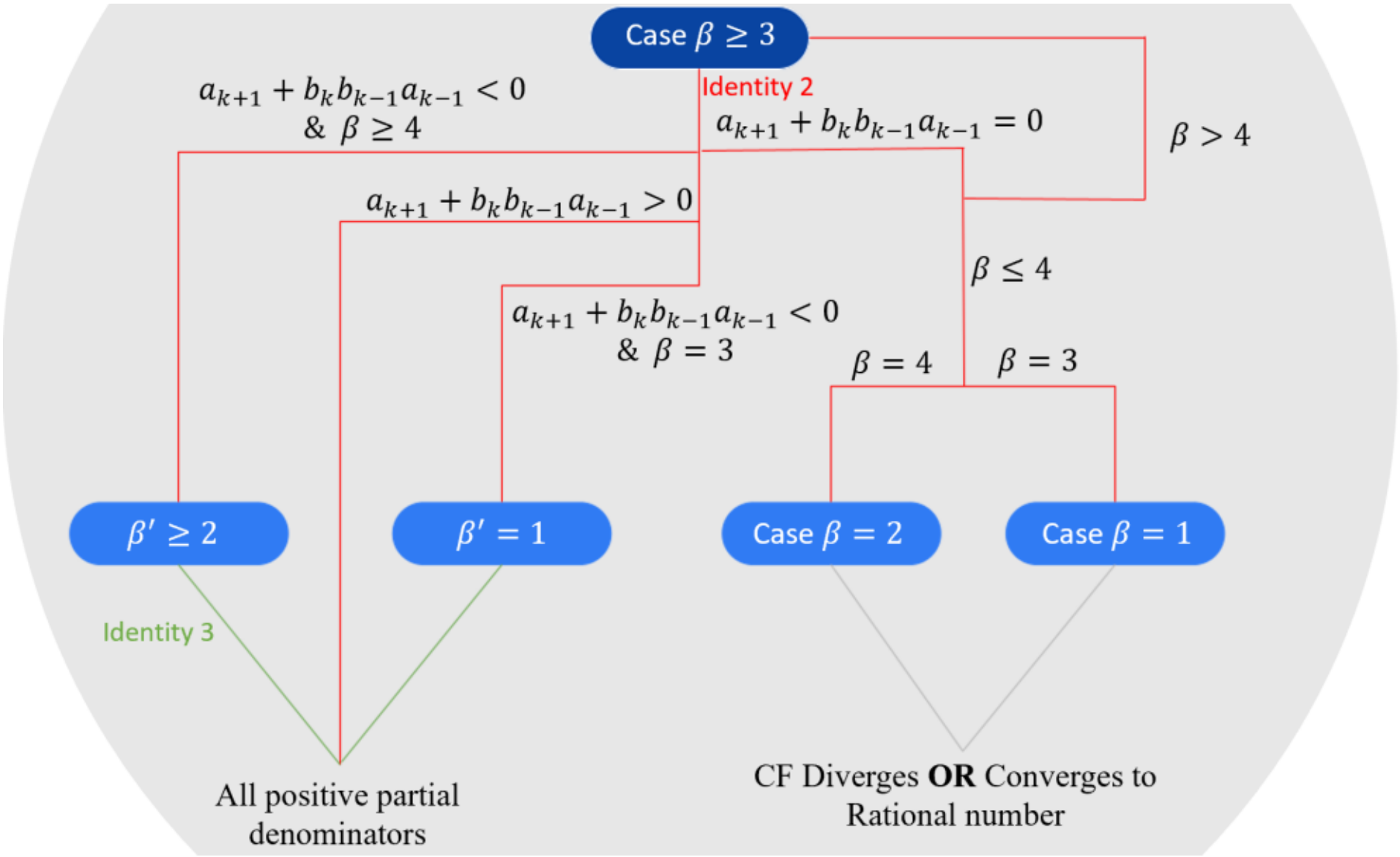

**Figure 5 | Flow chart proof of lemma 2**: Given a converging SICF with a single $a_k = a_{(n-1)\times\beta\,+\,k} = 0\ \forall n \in \mathbb{N}\backslash 0$ in each period, we show that there exists an equivalent representation (using identity 2 and 3) with all positive partial denominators. Recall, applying identity 2 decreases the period of the SICF by 2, therefore we are bound to either find a representation with all positive partial denominators or prove by contradiction that the SICF diverges or converges to a rational number.

We begin by assessing different periods.

<u>If $\beta = 1$:</u>

If $\beta = 1$, we obtain a diverging CF. We have a single layer meaning every layer has a 0 in its partial denominator and we obtain a CF, which satisfies:

$$\alpha = [(0, b_1), (0, b_1), \dots, (0, b_1), \dots] = \frac{b_1}{0+\frac{b_1}{0+\cdots}}.$$

As $b_1 = \pm 1$, at each $N^{th}$ approximation of the CF we obtain either $\infty$ or 0 and therefore the SICF diverges and is not relevant to our analysis.

<u>If $\beta = 2$:</u>

We obtain the following CF, which converges to a rational number 0,

$$\alpha = [(A_1(1), b_1), (0, b_2), \dots] = \lim_{N\to\infty}\left[\prod_{n=1}^{N}\begin{pmatrix}0 & b_1\\ 1 & A_1(n)\end{pmatrix}\begin{pmatrix}0 & b_2\\ 1 & 0\end{pmatrix}\right](0) = \lim_{N\to\infty}\left[\prod_{n=1}^{N}\begin{pmatrix}b_1 & 0\\ A_1(n) & b_2\end{pmatrix}\right](0),$$

Since $\begin{pmatrix}b_1 & 0\\ A_1(n) & b_2\end{pmatrix}(0) = 0$, if the limit exists it must be 0 which is rational and of no interest to us.

<u>If $\beta \geq 3$:</u>

As the period is $\geq 3$, we can "shift" the period such that the $a_k = 0$ element is between two matrices within the period and apply identity 2,

$$\alpha = [\dots, (a_{k-1}, b_{k-1}), (0, b_k), (a_{k+1}, b_{k+1}), \dots] \overset{\text{identity 2}}{=} [\dots, (a_{k+1} + b_k b_{k+1} a_{k-1},\ b_{k-1} b_k b_{k+1}), \dots]..$$

Where our period is reduced from $\beta$ to a new period $\beta' = \beta - 2$.
**If $\boldsymbol{a_{k+1} + b_k b_{k+1} a_{k-1} > 0}$,** we obtain a layer with a positive partial denominator and the total number of partial numerators equal to $-1$ either remains the same or decreases depending on the sign of $b_{k-1} b_k b_{k+1}$, satisfying lemma 2.

**If $\boldsymbol{a_{k+1} + b_k b_{k+1} a_{k-1} = 0}$,** then depending on our period we treat it by induction:

**If $\beta = 3$**, we obtain the case of $\beta = 1$ and therefore a CF which diverges.

**If $\beta = 4$**, we obtain the case of $\beta = 2$ and therefore a CF which converges to 0.

**If $\beta > 4$**, we obtain the case of $\beta \geq 3$ and treat it accordingly with the new period given by $\beta - 2$.

**If $a_{k+1} + b_k b_{k+1} a_{k-1} < 0$ and $\beta = 3$**, we obtain a SICF of period 1 (as the period is reduced to $\beta - 2$) with a negative partial denominator. This is referred to as $\beta' = 1$ in the flowchart. We can apply the equivalence theorem (identity 3 with $g_j = -1\ \forall j \in \mathbb{N}$ ) to ensure that we obtain a positive partial numerator:

$$\alpha = [\dots,(a_{k-1} - a_{k+1}, -b_{k-1}), \dots] = [\dots,(g_k(a_{k-1} - a_{k+1}), -g_k g_{k-1} b_{k-1}), \dots] = [\dots,(-(a_{k-1} - a_{j+1}), -b_{k-1}), \dots].$$

Again, we do not increase the number of $b_j's = -1$ as $g_k g_{k-1} = 1$.

**If $a_{k+1} + b_k b_{k+1} a_{k-1} < 0$ and $\beta \geq 4$** then we have another layer within the period and we could move the minus to, ensuring positive partial denominators. This is referred to as $\beta' \geq 2$ in the flow chart. This is done using the equivalence theorem (identity 3) with $g_{(n-1)\times\beta + k} = -1\ \forall n \in \{1,2,3 \dots\}$ and $g_j = 1\ \forall j \neq (n-1) \times \beta + k$.

$$\alpha = [\dots,(a_{k-1} - a_{k+1}, -b_{k-1}),(a_{k+2}, b_{k+2}), \dots] = [\dots,(g_k(a_{k-1} - a_{k+1}), -g_k b_{k-1}),(a_{k+2}, g_k b_{k+2}), \dots],$$

$$\alpha = [\dots,(-(a_{k-1} - a_{k+1}), b_{k-1}),(a_{k+2}, -b_{k+2}), \dots].$$

We do not increase the number the number of $b_j$'s $= -1$ in each period. If $b_{k+2} = -1$ we decrease the total amount of negative $b_j$'s, and if $b_{k+2} = 1$ we remain with the same amount of negative $b_j$'s as its given that $b_k b_{k+1} = -1$ and that was canceled out when we applied identity 3.

Thus, we ensured that all partial denominators are positive without increasing the number of $b_j$'s $= -1$ in the period. Furthermore, all identities were applied within a given period up to a shift which is equivalent to a Mobius transform, therefore the resultant CF with positive $a_j$'s will have an interlaced pattern.

Proof of Lemma 1

We show that, given a SICF that converges to an irrational number, if all $a_j$ are positive, we can always decrease the number of $b_j$'s$= -1$ in each period at least by 1. We first show, as base cases, that any SICF of period $\beta = 1$ or $\beta = 2$ can be represented as a simple interlaced CF up to some Mobius transform.

Base cases

The first base case is a SICF of period 1, which must be addressed as identity 1 requires at least two layers within a given period to be applied. The second base case is a SICF of period 2.

Base case: SICF of period $\beta = 1$

For any SICF of period 1, we shift the CF to some depth for which $A_1(n) > 2$ (for simplicity we keep $n = 1$ as starting index regardless of shift):

$$\alpha = \lim_{N\to\infty}\left[\prod_{n=1}^{N}\begin{pmatrix}0 & -1\\ 1 & A_1(n)\end{pmatrix}\right](0) = \lim_{N\to\infty}\left[\prod_{n=1}^{N}\begin{pmatrix}-1 & 0\\ 1 & 1\end{pmatrix}\begin{pmatrix}0 & 1\\ 1 & A_1(n)-1\end{pmatrix}\right](0) =,$$

$$= \lim_{N\to\infty}\left[\begin{pmatrix}-1 & 0\\ 1 & 1\end{pmatrix}\prod_{n=1}^{N}\begin{pmatrix}0 & 1\\ 1 & A_1(n)-1\end{pmatrix}\begin{pmatrix}-1 & 0\\ 1 & 1\end{pmatrix}\right](0) = \lim_{N\to\infty}\left[\begin{pmatrix}-1 & 0\\ 1 & 1\end{pmatrix}\prod_{n=1}^{N}\begin{pmatrix}0 & 1\\ 1 & A_1(n)-2\end{pmatrix}\begin{pmatrix}0 & 1\\ 1 & 1\end{pmatrix}\right](0).$$

We obtain a integer Mobius transform of a simple interlaced CF; Note, that we got a sub-sequence of the CF and as simple CFs always converge we obtain a simple interlaced expansion which converges to $\alpha$.

Base case: SICF of period $\beta = 2$

There are several possibilities of SICFs of period 2, we solve the possibility which seems most complex. The other cases are solved in a similar manner and often more simply. for example,

$$\alpha = \lim_{N\to\infty}\left[\prod_{n=1}^{N}\begin{pmatrix}0 & 1\\ 1 & A_1(n)\end{pmatrix}\begin{pmatrix}0 & -1\\ 1 & 1\end{pmatrix}\right](0) \overset{\text{identity 1}}{=}$$

$$\lim_{N\to\infty}\left[\prod_{n=1}^{N}\begin{pmatrix}0 & 1\\ 1 & A_1(n)-1\end{pmatrix}\begin{pmatrix}0 & 1\\ 1 & 1\end{pmatrix}\begin{pmatrix}0 & 1\\ 1 & 0\end{pmatrix}\right](0) \overset{\substack{\text{shift}\\ \text{identity 2}}}{=} \lim_{N\to\infty}\left[\begin{pmatrix}0 & 1\\ 1 & A_1(1)-1\end{pmatrix}\prod_{n=1}^{N}\begin{pmatrix}0 & 1\\ 1 & A_1(n+1)\end{pmatrix}\right](0).$$

A similar application of identity 1, a shift, then identity 2 can be easily applied in the case that $A_1(n)$'s layer has a negative partial numerator. We therefore analyze the more complicated case where all $b_j$'s$=-1$:

$$\alpha = \lim_{N\to\infty}\left[\prod_{n=1}^{N}\begin{pmatrix}0 & -1\\ 1 & A_1(n)\end{pmatrix}\begin{pmatrix}0 & -1\\ 1 & 1\end{pmatrix}\right](0).$$

We notice that $\begin{pmatrix}0 & -1\\ 1 & 1\end{pmatrix} = \begin{pmatrix}1 & -1\\ 0 & 1\end{pmatrix}\begin{pmatrix}1 & 0\\ 1 & 1\end{pmatrix}$, and therefore:

$$\alpha = \lim_{N\to\infty}\left[\prod_{n=1}^{N}\begin{pmatrix}0 & -1\\ 1 & A_1(n)\end{pmatrix}\begin{pmatrix}1 & -1\\ 0 & 1\end{pmatrix}\begin{pmatrix}1 & 0\\ 1 & 1\end{pmatrix}\right](0),$$

$$\alpha = \lim_{N\to\infty}\left[\begin{pmatrix}0 & -1\\ 1 & A_1(1)\end{pmatrix}\begin{pmatrix}1 & -1\\ 0 & 1\end{pmatrix}\prod_{n=1}^{N}\begin{pmatrix}1 & 0\\ 1 & 1\end{pmatrix}\begin{pmatrix}0 & -1\\ 1 & A_1(n+1)\end{pmatrix}\begin{pmatrix}1 & -1\\ 0 & 1\end{pmatrix}\right](0),$$

$$\alpha = \lim_{N\to\infty}\left[\begin{pmatrix}0 & -1\\ 1 & A_1(1)\end{pmatrix}\begin{pmatrix}1 & -1\\ 0 & 1\end{pmatrix}\prod_{n=1}^{N}\begin{pmatrix}0 & -1\\ 1 & A_1(n+1)-2\end{pmatrix}\right](0).$$

We obtain a 1-periodic SICF which can be represented as a simple interlaced CF (as shown in the base case of $\beta = 1$).

<u>General case: SICF of period β ≥ 3</u>

Given a SICF,

$$\alpha = [(a_1, b_1), (a_2, b_2), \ldots, (a_\beta, b_\beta), \ldots], \qquad b_j \in \{\pm 1\}, a_j > 0$$

There are 2 cases that must be addressed:

1. $\forall j \in \mathbb{N}\ b_j = -1$. Recall, we are discussing SICF with non-constant partial denominator sequences, therefore $\exists i \in \{1, \beta\}$ s.t $\deg(A_i) > 0$. Therefore, we can ensure $a_j > 1$ as $A_i(n) = a_{(n-1)\times\beta + k}$. As the period $\geq 3$ we can "shift" the period such that this $a_j$ is adjacent to another matrix from the left within the period. This enables us to apply identity 1 while guaranteeing that $a_j - 1 > 0$, This decreases the number of $b_j$'s$=-1$ in each period by 1.

$$\begin{pmatrix}0 & b_{j-1}\\ 1 & a_{j-1}\end{pmatrix}\begin{pmatrix}0 & -1\\ 1 & a_j\end{pmatrix} = \begin{pmatrix}0 & b_{j-1}\\ 1 & a_{j-1}-1\end{pmatrix}\begin{pmatrix}1 & 1\\ 0 & 1\end{pmatrix}\begin{pmatrix}-1 & 0\\ 1 & 1\end{pmatrix}\begin{pmatrix}0 & 1\\ 1 & a_j-1\end{pmatrix} = \begin{pmatrix}0 & b_{j-1}\\ 1 & a_{j-1}-1\end{pmatrix}\begin{pmatrix}0 & 1\\ 1 & 1\end{pmatrix}\begin{pmatrix}0 & 1\\ 1 & a_j-1\end{pmatrix}$$

We have therefore decreased the number of $b_j$'s$=-1$ and again treat potential zeros ($a_{j-1} - 1 = 0$) as outlined in lemma 2.

We can then apply the iterative decomposition as in case 2, as we guarantee that at least one partial numerator satisfies $b_j = 1$.

2. $\exists k \in \mathbb{N} | b_{k-1} = 1\ \&\ b_k = -1$. This case can occur in any SICF that has at least a single $b_j = 1$ in each period, as then we can shift the period to match the conditions order (first positive than negative layer) and then decompose the CF as using identity 1 and decrease the total number of $b_j$'s$=-1$ in each period.

$$\alpha = [\ldots, (a_{k-1}, 1), (a_k, -1), \ldots] \rightarrow \alpha = [\ldots, (a_{k-1}-1, 1), (1,1), (a_k - 1, 1), \ldots].$$

Note: this is done in every period of the SICF

We treat potential zeros on the left layer and right layer separately, enabling us to satisfy the assumption in Lemma 2; that there is only one zero partial denominator in each period.

**If $a_{k-1} - 1 = 0$**, we can simply add 1 to the previous layer, (guaranteed to have one as the period ≥ 3)

$$\alpha = [\ldots, (a_{k-2}+1, 1), (a_k - 1, 1), \ldots, (a_\beta, b_\beta), \ldots]$$

**If $\forall n \in \mathbb{N}\backslash 0, k = 2 + (n-1)\times\beta$** (meaning $a_{k-1}$ is the first layer of each period), we can shift the period and add 1 to the final layer.

$$\alpha = [\ldots, (a_k - 1, 1), \ldots, (a_\beta + 1, b_\beta), \ldots]$$

**If $a_k - 1 = 0$**, we utilize lemma 2 to return to a state in which all partial denominators are positive.

If we do not get any zero in the partial denominators, we increased the SICF's period from $\beta$ to $\beta+1$ and decreased the number of $b_j$'s$=-1$. We then continue to iteratively apply identity 1 from left to right in the period to every $b_j=-1$ we encounter. The iterative process continues until we either obtain a simple interlaced CF (reducing the number of $b_j$'s$=-1$ at each iteration until all $b_j$'s$=1$) or obtain a 0 that will be dealt with as formalized in the proof of Lemma 2.

Proof of convergence

Throughout the above proof of theorem 2 we utilize 3 identities of matrix multiplication which might change the limit of the CF. Identities 2 and 3 do not change the limit as identity 2 simply eliminates a sub-sequence of the CF and identity 3 does not affect convergence. Therefore, to prove the limit is unchanged we prove identity 1 does not affect the convergence of the SICF. The SICF converges therefore,

$$\alpha=[\dots,(a_{j-1}-1,1),(1,1),(a_j-1,1),\dots].$$

We define $N(n)=j+\beta n\ \forall n\in\mathbb{N}$ (as $N\to\infty$ means $n\to\infty$). The SICF satisfies:

$$\alpha=\lim_{N\to\infty}\frac{p_N}{q_N}=\lim_{N\to\infty}\begin{pmatrix}0&b_1\\1&a_1\end{pmatrix}\dots\begin{pmatrix}0&-1\\1&a_N\end{pmatrix}(0)=\lim_{N\to\infty}\begin{pmatrix}0&b_1\\1&a_1\end{pmatrix}\dots\begin{pmatrix}0&1\\1&a_{N-1}-1\end{pmatrix}\begin{pmatrix}0&1\\1&1\end{pmatrix}\begin{pmatrix}0&1\\1&a_N-1\end{pmatrix}(0).$$

We must assess whether this decomposition affects the convergence of the CF. For most subsequences we can trivially see there is no effect on the convergence, $\begin{pmatrix}0&1\\1&a_{N-1}-1\end{pmatrix}\begin{pmatrix}0&1\\1&1\end{pmatrix}(0)=\begin{pmatrix}0&1\\1&a_{N-1}\end{pmatrix}$. We therefore assess a non-trivial case:

$$\begin{pmatrix}0&b_1\\1&a_1\end{pmatrix}\dots\begin{pmatrix}0&1\\1&a_{N-1}-1\end{pmatrix}(0)=\begin{pmatrix}0&b_1\\1&a_1\end{pmatrix}\dots\begin{pmatrix}0&1\\1&a_{N-1}\end{pmatrix}(-1)$$

$$=\begin{pmatrix}p_{N-1}&p_N\\q_{N-1}&q_N\end{pmatrix}(-1)=\frac{p_N-p_{N-1}}{q_N-q_{N-1}}.$$

To prove that the above expression converges we must show that:

$$\lim_{N\to\infty}\frac{p_N-p_{N-1}}{q_N-q_{N-1}}=\lim_{N\to\infty}\frac{p_N}{q_N}=\alpha\to\lim_{N\to\infty}\frac{p_N-p_{N-1}}{q_N-q_{N-1}}-\frac{p_N}{q_N}=0.$$

We notice that the above expression satisfies:

$$\frac{p_N-p_{N-1}}{q_N-q_{N-1}}-\frac{p_N}{q_N}=\frac{q_N(p_N-p_{N-1})-p_N(q_N-q_{N-1})}{q_N(q_N-q_{N-1})}\overset{\det\left(\begin{pmatrix}p_{N-1}&p_N\\q_{N-1}&q_N\end{pmatrix}\right)=\pm1}{=}\frac{\pm1}{q_N(q_N-q_{N-1})}.$$

If the SICF converges to an irrational number, it guarantees that:

$$|q_N|\to_{N\to\infty}\infty$$

As if we assume that $q_N$ is bounded, meaning $|q_N|<B\ \ B\in\mathbb{Z}$, then $p_N$ must be bounded and then we would have a sub-sequence that converges to a rational number.

However, we might encounter problems for cases which $q_N=q_{N-1}$ for an infinite amount of $N$'s. Since, $\det\left(\begin{pmatrix}p_{N-1}&p_N\\q_{N-1}&q_N\end{pmatrix}\right)=\pm1$, if $q_N=q_{N-1}$ then they must be $\pm1$ for infinity many $N$'s which contradicts $|q_N|\to_{N\to\infty}\infty$. Therefore, it is not possible for $q_N=q_{N-1}$ for an infinite amount of $\mathrm{N}'$s and we obtain:

$$\lim_{k\to\infty}|q_N(q_N-q_{N-1})|=\infty\to\lim_{k\to\infty}\frac{\pm1}{q_N(q_N-q_{N-1})}=0.$$

Therefore,

$$\lim_{N\to\infty}\frac{p_N-p_{N-1}}{q_N-q_{N-1}}=\lim_{N\to\infty}\frac{p_N}{q_N}=\alpha$$

And the resultant CF converges to the same limit.

Overall, lemma 1 guarantees that given a SICF (all $a_j$ are positive) it can always reduce the number of $b_j$'s$=-1$ in each period yet it may result in edge cases of 0 in the partial denominator. While lemma 2 guarantees if all $a_j$ in each period are positive but one, $\exists k\in\mathbb{N}|a_k=a_{k+(n-1)\times\beta}=0$, there is an equivalent representation with all $a_j$ being positive without increasing the number of $b_j$'s$=-1$ in each period. Both the lemma's show that we can ensure that the number of $b_j$'s$= -1$ is monotonically decreasing ensuring that for a converging SICF, we can find a representation such that all $b_j$'s$= 1$ while all $a_j$'s remain positive. Therefore, by proving both the induction hypothesis we prove that every constant equal to a SICF has a Mobius transform equal to a simple CF.

## Appendix C.2: Example of moving from an SICF to a Simple Interlaced CF

We take the following conjecture of $\frac{2}{\tan(1)}$ found by the ESMA algorithm,

$$\frac{2}{\tan(1)} - 2 = \lim_{N\to\infty}\left[\prod_{n=1}^{N}\begin{pmatrix}0 & -1\\ 1 & 3n-1\end{pmatrix}\begin{pmatrix}0 & -1\\ 1 & 2\end{pmatrix}\begin{pmatrix}0 & -1\\ 1 & 3n\end{pmatrix}\begin{pmatrix}0 & -1\\ 1 & 2+12n\end{pmatrix}\right](0). \quad (34)$$

We attempt to transform it to a simple interlaced CF, using the method from the constructive proof described in Appendix C.1:

$$\frac{2}{\tan(1)} - 2 = \lim_{N\to\infty}[\prod_{n=1}^{N}\begin{pmatrix}-1 & 0\\ 1 & 1\end{pmatrix}\begin{pmatrix}0 & 1\\ 1 & 3n-2\end{pmatrix}\begin{pmatrix}0 & -1\\ 1 & 2\end{pmatrix}\begin{pmatrix}0 & -1\\ 1 & 3n\end{pmatrix}\begin{pmatrix}0 & -1\\ 1 & 1+12n\end{pmatrix}\begin{pmatrix}1 & 1\\ 0 & 1\end{pmatrix}](0),$$

$$= \lim_{N\to\infty}[\begin{pmatrix}-1 & 0\\ 1 & 1\end{pmatrix}\prod_{n=1}^{N}\begin{pmatrix}0 & 1\\ 1 & 3n-3\end{pmatrix}\begin{pmatrix}1 & 1\\ 0 & 1\end{pmatrix}\begin{pmatrix}-1 & 0\\ 1 & 1\end{pmatrix}\begin{pmatrix}0 & 1\\ 1 & 1\end{pmatrix}\begin{pmatrix}0 & -1\\ 1 & 3n\end{pmatrix}\begin{pmatrix}0 & -1\\ 1 & 1+12n\end{pmatrix}\begin{pmatrix}0 & 1\\ 1 & 1\end{pmatrix}](0),$$

$$= \lim_{N\to\infty}[\begin{pmatrix}-1 & 0\\ 1 & 1\end{pmatrix}\prod_{n=1}^{N}\begin{pmatrix}0 & 1\\ 1 & 3n-3\end{pmatrix}\begin{pmatrix}0 & 1\\ 1 & 1\end{pmatrix}\begin{pmatrix}0 & 1\\ 1 & 1\end{pmatrix}\begin{pmatrix}0 & -1\\ 1 & 3n\end{pmatrix}\begin{pmatrix}0 & -1\\ 1 & 1+12n\end{pmatrix}\begin{pmatrix}0 & 1\\ 1 & 1\end{pmatrix}](0),$$

$$= \lim_{N\to\infty}[\begin{pmatrix}-1 & 0\\ 1 & 1\end{pmatrix}\prod_{n=1}^{N}\begin{pmatrix}0 & 1\\ 1 & 3n-3\end{pmatrix}\begin{pmatrix}0 & 1\\ 1 & 1\end{pmatrix}\begin{pmatrix}0 & 1\\ 1 & 0\end{pmatrix}\begin{pmatrix}0 & 1\\ 1 & 1\end{pmatrix}\begin{pmatrix}0 & 1\\ 1 & 3n-1\end{pmatrix}\begin{pmatrix}0 & -1\\ 1 & 1+12n\end{pmatrix}\begin{pmatrix}0 & 1\\ 1 & 1\end{pmatrix}](0),$$

$$= \lim_{N\to\infty}[\begin{pmatrix}-1 & 0\\ 1 & 1\end{pmatrix}\prod_{n=1}^{N}\begin{pmatrix}0 & 1\\ 1 & 3n-3\end{pmatrix}\begin{pmatrix}0 & 1\\ 1 & 1\end{pmatrix}\begin{pmatrix}0 & 1\\ 1 & 0\end{pmatrix}\begin{pmatrix}0 & 1\\ 1 & 1\end{pmatrix}\begin{pmatrix}0 & 1\\ 1 & 3n-2\end{pmatrix}\begin{pmatrix}0 & 1\\ 1 & 1\end{pmatrix}\begin{pmatrix}0 & 1\\ 1 & 12n\end{pmatrix}\begin{pmatrix}0 & 1\\ 1 & 1\end{pmatrix}](0),$$

$$= \lim_{N\to\infty}[\begin{pmatrix}-1 & 0\\ 1 & 1\end{pmatrix}\prod_{n=1}^{N}\begin{pmatrix}0 & 1\\ 1 & 3n-3\end{pmatrix}\begin{pmatrix}0 & 1\\ 1 & 2\end{pmatrix}\begin{pmatrix}0 & 1\\ 1 & 3n-2\end{pmatrix}\begin{pmatrix}0 & 1\\ 1 & 1\end{pmatrix}\begin{pmatrix}0 & 1\\ 1 & 12n\end{pmatrix}\begin{pmatrix}0 & 1\\ 1 & 1\end{pmatrix}](0),$$

$$= \lim_{N\to\infty}[\begin{pmatrix}-1 & 0\\ 1 & 1\end{pmatrix}\begin{pmatrix}0 & 1\\ 1 & 0\end{pmatrix}\prod_{n=1}^{N}\begin{pmatrix}0 & 1\\ 1 & 2\end{pmatrix}\begin{pmatrix}0 & 1\\ 1 & 3n-2\end{pmatrix}\begin{pmatrix}0 & 1\\ 1 & 1\end{pmatrix}\begin{pmatrix}0 & 1\\ 1 & 12n\end{pmatrix}\begin{pmatrix}0 & 1\\ 1 & 1\end{pmatrix}\begin{pmatrix}0 & 1\\ 1 & 3n\end{pmatrix}](0),$$

$$-\frac{1}{\frac{2}{\tan(1)}-2} - 1 = \lim_{N\to\infty}[\prod_{n=1}^{N}\begin{pmatrix}0 & 1\\ 1 & 2\end{pmatrix}\begin{pmatrix}0 & 1\\ 1 & 3n-2\end{pmatrix}\begin{pmatrix}0 & 1\\ 1 & 1\end{pmatrix}\begin{pmatrix}0 & 1\\ 1 & 12n\end{pmatrix}\begin{pmatrix}0 & 1\\ 1 & 1\end{pmatrix}\begin{pmatrix}0 & 1\\ 1 & 3n\end{pmatrix}](0).$$

## Appendix D.1: Proof of Theorem 3 Predicting Degrees of Polynomial CFs semi-equivalent to simple Interlaced CFs

**Theorem 3.** Given a simple interlaced continued fraction satisfying $B_i(n) = 1\ \&\ A_i(n) > 0$, where $A_i, B_i \in \mathbb{Z}[x]$, $i \in \{1,..,\beta\}, \forall n \in \mathbb{N}$ of period $\beta$, where $\exists i \in \{1,..,\beta\}$ s.t $\deg(A_i) > 0$, its semi-equivalent polynomial continued fraction's partial numerator $b_n'$ and partial denominator $a_n'$ ($\forall n \in \mathbb{N}$) have the following polynomial degrees:

$$\deg(b') = \sum_{i=1}^{\beta-1} 2\times\deg(A_i),\ \deg(a') = \left[\sum_{i=1}^{\beta-1} 2\times\deg(A_i)\right] + \deg(A_\beta).$$

Complementary theorem

As defined previously, the form of a simple interlaced CF, $\alpha$, with a period of $\beta$ is given by:

$$\alpha = \lim_{N\to\infty}\left[\prod_{n=1}^{N}\prod_{i=1}^{\beta}\begin{pmatrix}0 & B_i(n)\\ 1 & A_i(n)\end{pmatrix}\right](0) = \lim_{N\to\infty}[\prod_{n=1}^{N}M_n](0) = \lim_{N\to\infty}\left[\prod_{n=1}^{N}\begin{pmatrix}c_n & d_n\\ e_n & f_n\end{pmatrix}\right](0),$$

$$A_i \in \mathbb{Z}[x],\ A_i(n) > 0,\ \forall n \in \mathbb{N},\ i \in \{1,\dots,\beta\}.$$

We prove a complementary theorem which will enable the proof of theorem 3.

**Complementary Theorem.** The collapsed matrix sequences of a simple interlaced CF of period $\beta \geq 2$, have the following polynomial degree with $n$:

$$(***)\begin{cases} \deg(c) = \begin{cases} \sum_{i=2}^{\beta-1} \deg(A_i) & if\ \beta > 2 \\ 0 & \end{cases} \\ \deg(d) = \sum_{i=2}^{\beta} \deg(A_i) \\ \deg(e) = \sum_{i=1}^{\beta-1} \deg(A_i) \\ \deg(f) = \sum_{i=1}^{\beta} \deg(A_i) \end{cases} \quad (35)$$

**Proof:**

We prove the complementary theorem using induction.

<u>Base case of induction $\beta = 2$:</u>

For $\beta = 2$, we obtain the following collapsed matrix:

$$\begin{pmatrix} 0 & 1 \\ 1 & A_1(n) \end{pmatrix}\begin{pmatrix} 0 & 1 \\ 1 & A_2(n) \end{pmatrix} = \begin{pmatrix} 1 & A_2(n) \\ A_1(n) & A_1(n)A_2(n) + 1 \end{pmatrix} = \begin{pmatrix} c_n & d_n \\ e_n & f_n \end{pmatrix}.$$

Meaning the induction hypothesis is met.

We assume the induction hypothesis is met for all simple interlaced CFs of period $\leq \beta - 1$ and prove that it is met for all simple interlaced CF of period $\beta$. We consider the effect of multiplying some product of the first $\beta - 1$ layer to the $\beta^{\text{th}}$ layer of the $n^{\text{th}}$ collapsed matrix:

$$\begin{pmatrix} C_n & D_n \\ E_n & F_n \end{pmatrix}\begin{pmatrix} 0 & 1 \\ 1 & A_\beta(n) \end{pmatrix} = \begin{pmatrix} D_n & C_n + D_n A_\beta(n) \\ F_n & E_n + F_n A_\beta(n) \end{pmatrix} = \begin{pmatrix} c_n & d_n \\ e_n & f_n \end{pmatrix}$$

By our assumption, the product of the first $\beta - 1$ layers have the following polynomial degrees:

$$\deg(C) = \sum_{i=2}^{\beta-2} \deg(A_i),\ \deg(D) = \sum_{i=2}^{\beta-1} \deg(A_i),\ \deg(E) = \sum_{i=1}^{\beta-2} \deg(A_i),\ \deg(F) = \sum_{i=1}^{\beta-1} \deg(A_i).$$

Notice that as the elements of the collapsed matrix are a product and linear combination of positive polynomial sequences the elements are also positive polynomial sequences: $c, d, e, f \in \mathbb{Z}[x],\ c_n \geq 0, d_n, e_n, f_n > 0$ and therefore there can be no cancelation of polynomials. The resultant collapsed matrix satisfies:

$$\deg(c) = \sum_{i=2}^{\beta-1} \deg(A_i),\ \deg(d) = \sum_{i=2}^{\beta} \deg(A_i),\ \deg(e) = \sum_{i=1}^{\beta-1} \deg(A_i), \text{and } \deg(f) = \sum_{i=1}^{\beta} \deg(A_i)$$

Proving our complementary theorem.

**<u>Proof of Theorem 3 (and lemma 1)</u>**

We utilize the complementary theorem to analyze the resultant polynomial CFs obtained when applying the Folding transform on the $n^{\text{th}}$ collapsed matrix. Given a collapsed matrix ($n > 1$) for some interlaced CF, the Folding transform on the $n^{\text{th}}$ collapsed matrix is given by:

$$\begin{pmatrix} e_n & -c_n \\ 0 & 1 \end{pmatrix}\begin{pmatrix} c_n & d_n \\ e_n & f_n \end{pmatrix}\begin{pmatrix} 1 & c_{n+1} \\ 0 & e_{n+1} \end{pmatrix} = \begin{pmatrix} 0 & e_{n+1}(e_n d_n - c_n f_n) \\ e_n & e_n c_{n+1} + f_n e_{n+1} \end{pmatrix} =_{\text{equivalence transform}}$$

$$\begin{pmatrix} 0 & e_{n-1}e_{n+1}(e_n d_n - c_n f_n) \\ 1 & (e_n c_{n+1} + f_n e_{n+1}) \end{pmatrix}.$$

Utilizing the determinant property (see Appendix B.4):

$$\left|\begin{pmatrix} c_n & d_n \\ e_n & f_n \end{pmatrix}\right| = (-1)^\beta \prod_{i=1}^{\beta} B_i.$$

In this case $B_i = 1\ \forall i \in \mathbb{N}$, we therefore obtain:

$$-\left|\begin{pmatrix} c_n & d_n \\ e_n & f_n \end{pmatrix}\right| = e_n d_n - c_n f_n = (-1)^{\beta+1} = \begin{cases} 1 & \beta(\text{mod})2 = 1 \\ -1 & \beta(\text{mod})2 = 0 \end{cases}\text{'}$$

And therefore, we obtain:

$$\begin{pmatrix} 0 & e_{n-1}e_{n+1}(e_n d_n - c_n f_n) \\ 1 & (e_n c_{n+1} + f_n e_{n+1}) \end{pmatrix} = \begin{pmatrix} 0 & e_{n-1}e_{n+1}(-1)^{\beta+1} \\ 1 & (e_n c_{n+1} + f_n e_{n+1}) \end{pmatrix}.$$

As $c_n \geq 0, d_n, e_n, f_n > 0 \; \forall n > 1$ we can deduce that $(e_n c_{n+1} + f_n e_{n+1}) > 0$ thus **proving lemma 3**. Therefore, the resultant polynomial degree is the sum of polynomial degrees of the polynomial products (recalling equation (26) in Appendix B.2):

$$(M_1 U_2)^{-1}(\alpha) = \lim_{N\to\infty} \left[ \prod_{n=2}^{N} \begin{pmatrix} 0 & e_{n-1}e_{n+1}(-1)^{\beta+1} \\ 1 & (e_n c_{n+1} + f_n e_{n+1}) \end{pmatrix} \right] (0) \, .$$

If we denote $b'_n$ and $a'_n$ as the partial numerator and denominator of $\alpha'$, respectively,

$$\deg(b') = 2\deg(e) = \textstyle\sum_{i=1}^{\beta-1} 2 * \deg(A_i), \; \deg(a') = \deg(e) + \deg(f) = \sum_{i=1}^{\beta-1} 2 * \deg(A_i) + \deg(A_\beta). \qquad (36)$$

Note, $\deg(e) + \deg(f) \geq \deg(e) + \deg(c)$ as $\deg(f) = \deg(c) + \deg(A_1) + \deg(A_\beta)$ by $(***)$.

Thus, **proving Theorem 3.**

## Appendix D.2: Proof of Corollary 2: A Polynomial CF semi-equivalent to a simple Interlaced Continued Fraction converges super-exponentially

**Corollary 2.** A polynomial continued fraction semi-equivalent to a simple interlaced continued fraction, satisfying $\exists i \in \{1, .., \beta\}$ s.t $\deg(A_i) > 0$, converges super-exponentially.

**Proof:**

In paper [6] the following condition for super-exponential convergence of a polynomial CF was proven: Given a polynomial CF with partial numerator $b'_n$ and partial denominator $a'_n$ the polynomial CF converges super-exponentially iff $\frac{\deg(b')}{\deg(a')} < 2$. In the case of a simple interlaced CF, its semi-equivalent polynomial CF satisfies:

$$\frac{\deg(b')}{\deg(a')} = \frac{\sum_{i=1}^{\beta-1} 2\times\deg(A_i)}{\sum_{i=1}^{\beta-1} 2\times\deg(A_i) + \deg(A_\beta)} \leq \frac{\sum_{i=1}^{\beta-1} 2\times\deg(A_i)}{\sum_{i=1}^{\beta-1} 2\times\deg(A_i)} = 1 < 2.$$

Therefore, the polynomial CF converges at a super-exponential rate.

## Appendix D.3: Proof of Corollary 3: Irrationality of a Simple Interlaced Continued Fraction

**Corollary 3.** Any number $\alpha \in \mathbb{R}$ equal to a simple interlaced continued fraction converges to an irrational limit.

**Proof:**

Utilizing Tietze's Criterion it can be directly proven that any simple interlaced CF is irrational. We denote α as the constant equal to the simple interlaced CF with partial numerator $b_n$ and partial denominator $a_n$, and α′ as its semi-equivalent constant equal to the polynomial CF. Recall Tietze's Criterion (as outlined in the second page of [20]):

Let $\{b_j\}_{j=1}^{\infty}$ be a sequence of integers and $\{a_j\}_{j=1}^{\infty}$ a sequence of positive integers, if there exists a positive integer $N_0$ such that

$$\begin{cases} a_j \geq |b_j| \\ a_j \geq |b_j| + 1 \;\; \text{for } b_{j+1} < 0 \end{cases} \quad \forall j \geq N_0. \qquad (37)$$

Then $\alpha' = \frac{b_1}{a_1 + \frac{b_2}{a_2 + \cdots}}$ converges and its limit is irrational.

By definition, a simple interlaced CF satisfies $a_j \geq |b_j| = 1$ (by lemma 1) and $b_j > 0 \; \forall j$ and therefore a simple Interlaced CF converges to an irrational limit.

## Appendix E: Interesting Continued Fractions found with the Folding Transform

Representation of the Golden Ratio as a Balanced Polynomial CF

Looking at the following SICF, where $k, q, r \in \mathbb{Z}$:

$$\alpha = \lim_{N\to\infty}\left[\prod_{n=1}^{N}\begin{pmatrix}0 & -1\\ 1 & 1\end{pmatrix}\begin{pmatrix}0 & 1\\ 1 & kn+q\end{pmatrix}\begin{pmatrix}0 & -1\\ 1 & kn+r\end{pmatrix}\begin{pmatrix}0 & -1\\ 1 & 1\end{pmatrix}\right](0) =$$
$$\lim_{N\to\infty}\left[\begin{pmatrix}0 & -1\\ 1 & 1\end{pmatrix}\begin{pmatrix}0 & 1\\ 1 & k+q\end{pmatrix}\prod_{n=1}^{N}\begin{pmatrix}0 & 1\\ 1 & 1+k+q-r\end{pmatrix}\right](0).$$

This converges for $(1+k+q-r)^2 > -4$ [6], we constrain the above expression to obtain the CF of Golden Ratio $r = 1\ k = 1\ q = 0$:

$$\varphi - 1 = \lim_{N\to\infty}\left[\prod_{n=1}^{N}\begin{pmatrix}0 & 1\\ 1 & 1\end{pmatrix}\right](0)$$

Meaning that looking at the original SICF:

$$-\frac{1}{\varphi} = \begin{pmatrix}0 & -1\\ 1 & 1\end{pmatrix}(\varphi - 1) = \lim_{N\to\infty}\left[\begin{pmatrix}0 & -1\\ 1 & 1\end{pmatrix}\begin{pmatrix}0 & 1\\ 1 & 1\end{pmatrix}\prod_{n=1}^{N}\begin{pmatrix}0 & 1\\ 1 & 1\end{pmatrix}\right](0),$$

However, we can also apply the Folding transform to the SICF before collapsing it to a 1-periodic CF, $r = 1, k = 1, q = 0$ so we obtain:

$$-\frac{1}{\varphi} = \lim_{N\to\infty}\left[\prod_{n=1}^{N}\begin{pmatrix}0 & -1\\ 1 & 1\end{pmatrix}\begin{pmatrix}0 & 1\\ 1 & n\end{pmatrix}\begin{pmatrix}0 & -1\\ 1 & n+1\end{pmatrix}\begin{pmatrix}0 & -1\\ 1 & 1\end{pmatrix}\right](0).$$

We can collapse the above layers and apply the Folding transform to obtain the following polynomial CF:

$$\begin{pmatrix}21 & 15\\ 3 & 1\end{pmatrix}\left(-\frac{1}{\varphi}\right) = \lim_{N\to\infty}\left[\prod_{n=2}^{N}\begin{pmatrix}0 & n^4+4n^3+2n^2-4n-3\\ 1 & -n^2-3n-3\end{pmatrix}\right](0).$$

As a result, we found a balanced polynomial CF representation for the Golden Ratio:

$$\frac{21\left(-\frac{1}{\varphi}\right)+15}{3\left(-\frac{1}{\varphi}\right)+1} = \frac{-21+15\varphi}{-3+\varphi} = \lim_{N\to\infty}\left[\prod_{n=2}^{N}\begin{pmatrix}0 & n^4+4n^3+2n^2-4n-3\\ 1 & -n^2-3n-3\end{pmatrix}\right](0).$$

By rationalizing the denominator, we obtain a simpler expression:

$$\frac{-21+15\varphi}{-3+\varphi} = \frac{60-48\sqrt{5}}{-30} = \frac{-30+24\sqrt{5}}{15} = \frac{-54+\frac{48(1+\sqrt{5})}{2}}{15} = \frac{-54+48\varphi}{15}$$

$$\frac{-54+48\varphi}{15} = \lim_{N\to\infty}\left[\prod_{n=2}^{N}\begin{pmatrix}0 & n^4+4n^3+2n^2-4n-3\\ 1 & -n^2-3n-3\end{pmatrix}\right](0). \quad (38)$$